\documentclass{article}

\pdfoutput=1

\usepackage{mysetup}
\usepackage{myarticle}
\usepackage{mysymbols}
\usepackage{mycustomsetup}

\graphicspath{{figures/}}

\usepackage[
    backref=true,
    backend=biber,
    giveninits=true,
    style=authoryear,
    uniquename=init,
    maxbibnames=99,
    date=year
]{biblatex}
\begin{document}

\title{Combinatorial Optimization enriched Machine Learning to solve the Dynamic Vehicle Routing Problem with Time Windows}

\newcommand\CoAuthorMark{\footnotemark[\arabic{footnote}]} 
\author[1]{Léo Baty\footnote{The first three authors contributed equally to this work.}}
\author[2]{Kai Jungel\protect\CoAuthorMark}
\author[2]{Patrick S. Klein\protect\CoAuthorMark}
\author[1]{Axel Parmentier}
\author[2,3]{Maximilian Schiffer}

\affil[1]{\small CERMICS, Ecole des Ponts, Marne-la-Vallée, France}
\affil[2]{\small School of Management, Technical University of Munich, Munich, Germany}
\affil[3]{\small Munich Data Science Institute, Technical University of Munich, Munich, Germany}

\date{}

\maketitle

\begin{abstract}
    With the rise of e-commerce and increasing customer requirements, logistics service providers face a new complexity in their daily planning, mainly due to efficiently handling same day deliveries. Existing multi-stage stochastic optimization approaches that allow to solve the underlying dynamic vehicle routing problem are either computationally too expensive for an application in online settings, or -- in the case of reinforcement learning -- struggle to perform well on high-dimensional combinatorial problems. To mitigate these drawbacks, we propose a novel machine learning pipeline that incorporates a combinatorial optimization layer. We apply this general pipeline to a dynamic vehicle routing problem with dispatching waves, which was recently promoted in the EURO Meets NeurIPS Vehicle Routing Competition at NeurIPS 2022. Our methodology ranked first in this competition, outperforming all other approaches in solving the proposed dynamic vehicle routing problem. With this work, we provide a comprehensive numerical study that further highlights the efficacy and benefits of the proposed pipeline beyond the results achieved in the competition, e.g., by showcasing the robustness of the encoded policy against unseen instances and scenarios.

\end{abstract}

\clearpage

\section{Introduction}\label{sec:introduction}

With the rise of e-commerce during the last decade, \glspl*{LSP} were exposed to increasing customer requirements, particularly with respect to (fast) delivery times. Accordingly, the concept of same-day deliveries, where \glspl*{LSP} guarantee to fulfill an order on the day on which they receive it, became a key element in the B2C sector. In fact, the share of same-day deliveries grew by 18.5\% in 2021~\parencite{prnewswireSddGrowth}. As e-commerce and B2C deliveries remain competitive markets with several major players, retailers and \glspl*{LSP} continuously aim to outbid each other, which led to continuously shortened lead times, offering delivery within as little as two hours for certain product types in selected cities~\parencite{twoHourSdd}. 

Realizing last-mile deliveries within such short planning horizons remains inherently challenging from an efficiency perspective as \glspl*{LSP} generally trade off short lead times against oversized resources, e.g., by maintaining an oversized fleet to always be able to immediately react to incoming orders. In fact, the concept of same-day deliveries leads the concept of day-ahead planning, which has been the status quo in last-mile logistics for decades, ad absurdum. Instead of solving a combinatorially complex but static planning problem to determine cost-efficient delivery routes, \glspl*{LSP} have to dynamically dispatch orders to vehicles and route these vehicles in an online problem setting. Taking decisions in such a setting requires anticipating the benefit of dispatching or delaying an order while still inheriting the combinatorial complexity of the corresponding static planning problem and handling uncertainty with respect to future incoming orders. 

The challenges that arise in such planning problems in practice invigorate the interest in dynamic~\gls*{VRP} variants from a scientific perspective. State-of-the-art methodologies to solve such problems model the underlying planning task as a multi-stage stochastic optimization problem solved either directly~\parencite{pillac2013, Soeffker2022} or via reinforcement learning~\parencite{NazariReinforcementlearning2018, hildebrandtThomasUlmer2023, bassoKulcsarDiazQu2022}. However, both of these approaches bear a major drawback. Reinforcement learning based algorithms succeed in taking anticipating decisions but often struggle when being applied to high-dimensional combinatorial problems. Contrarily, stochastic optimization techniques are generally amenable to combinatorial problem settings but struggle with respect to computational efficiency in high-dimensional problems, which makes them impracticable to use in a dynamic setting~\parencite{PflugMultistage2014,CarpentierMultistage2015}.

Against this background, we propose a new methodological approach that mitigates the aforementioned shortcomings. Specifically, we develop a \gls*{ML} pipeline with an integrated \gls*{CO}-layer that allows to efficiently solve the dynamic~\gls*{VRP}. 
This pipeline mitigates the challenges of multi-stage stochastic optimization problems by design: its \gls*{ML}-layer allows to incorporate uncertainty by adequately parameterizing an instance of the underlying deterministic \gls*{CO} problem, which can then be efficiently solved within the \gls*{CO}-layer. We used this pipeline in the EURO Meets NeurIPS Vehicle Routing Competition at NeurIPS 2022 \parencite{euromeetsneurips2022}, where it outperformed all other approaches.

\subsection{Related work}\label{subsec:literature-review}
Our work contributes to two different streams of research: from an application perspective it relates to dynamic~\glspl*{VRP} and from a methodological perspective it relates to \gls*{CO}-enriched \gls*{ML}. We provide an overview of both related research streams in the following. For a general overview of \glspl*{VRP} we refer to \cite{Vidal2020} and for a general overview of CO-enriched ML we refer to \cite{bengioMachine2021}, and \cite{kotaryEndEnd2021}.

\paragraph{Dynamic Vehicle Routing Problems}
~Dynamic \glspl*{VRP} account for the dynamic nature of real-world processes where some problem data, such as the customers to serve, their demand, or the travel time between them, is not known in advance but revealed over time. Dynamic~\glspl*{VRP} hence have a diverse field of applications ranging from ride-hailing \parencite{Jungel2023} over grocery delivery \parencite{fikar2018} to emergency services \parencite{ALINAGHIAN2019}. To keep this literature review concise, we focus on dynamic~\glspl*{VRP} in the context of dynamic dispatching problems in the following and refer to~\textcite{pillac2013, humberto2021}~, and~\textcite{Ulmer2018} for comprehensive reviews of the field. Approaches to solve dynamic~\glspl*{VRP} can be broadly categorized into \gls*{CO}-based approaches, which leverage the combinatorial structure of the underlying problem, and \gls*{ML}-based approaches, which learn prescient policies accounting for the uncertainty of future customers.

Pure \gls*{CO} approaches generally amend solution methods developed for static \glspl*{VRP} to the dynamic case. Specifically, these \gls*{CO} approaches embed (meta-)heuristics into rolling horizon frameworks, i.e., solve a static variant of the considered problem each time new information enters the system~\parencite[see, e.g.,][]{ritzinger2016, humberto2021}. Here, \emph{myopic} approaches \parencite[see, e.g.,][]{gendreau1999, steever2019} utilize only the current problem state, while \emph{look-ahead} approaches \parencite[e.g.,][]{Flatberg2007, bent2004} take into account potential realizations of future periods, e.g., via sampling.

\gls*{ML} approaches learn policies which account for the uncertainty of future observations. Accordingly, they model the underlying dynamic problem as a markov decision process solved with either approximate dynamic programming methods~\parencite[see, e.g.,][]{Ulmer2018}, i.e., policy- or value function approximation, or reinforcement learning~\parencite[see, e.g.,][]{NazariReinforcementlearning2018, kool2018Attention, JoeLearning2020}. We refer to~\cite{RazaLearningbook2022} and \cite{hildebrandtThomasUlmer2023} for a review on reinforcement learning applied to dynamic~\glspl*{VRP}.

As can be seen, various works exist that solve dynamic~\glspl*{VRP}. Here, most approaches either utilize classical \gls*{CO} algorithms by sampling future scenarios or apply \gls*{ML} to approximate decision values which account for future expected rewards. All of these approaches contain at minimum one of the following shortcomings: classical \gls*{CO}-based algorithms struggle to amend to real-time requirements of dynamic problem settings, while \gls*{ML}-based approaches often lack solution quality as they do not take the problem's combinatorial structure into account.
A truly integrated approach which combines \gls*{CO} and \gls*{ML}, thus leveraging the advantage of \gls*{ML} in dynamic settings without disregarding the problem's combinatorial structure, has so far not been proposed for dynamic \glspl*{VRP}.

\paragraph{Combinatorial Optimization enriched Machine Learning}
Many real-world combinatorial dispatching problems are subject to uncertain future events.
To find combinatorial solutions which account for these uncertain events, one can integrate \gls*{CO}-layers into \gls*{ML}-based pipelines. We refer to such pipelines as \gls*{CO}-enriched \gls*{ML} pipelines. The main obstacle for using \gls*{CO}-layers in \gls*{ML} pipelines is their piecewise constant nature.
Specifically, gradients are zero almost everywhere, and thus uninformative in such settings, rendering straightforward backpropagation ineffective.
State-of-the-art methods address this issue and introduce regularization techniques that smoothen \gls*{CO}-layers to enable meaningful gradient computation, allowing their usage in \gls*{ML}-based pipelines.
Both additive~\parencite{berthetLearning2020} and multiplicative~\parencite{dalleLearning2022} regularization approaches have been applied successfully to \gls*{CO}-enriched \gls*{ML} pipelines in supervised learning settings with Fenchel-Young losses~\parencite{blondelLearning2020}.

A common application of these pipelines are hard \gls*{CO} problems, e.g., single machine scheduling with release dates~\parencite{parmentierLearning2021-1}, for which classical \gls*{CO} approaches are often intractable.
Here, the statistical model learns to parameterize an embedded, tractable \gls*{CO} problem that shares the feasible solution space with the intractable hard problem, such that it produces a solution that is valid for the original problem. This approach can be adapted to learn the parameterization of multi-stage optimization problems such as the two-stage stochastic minimum spanning tree problem~\parencite{dalleLearning2022} or the stochastic vehicle scheduling problem~\parencite{parmentierLearning2021, parmentierLearningApproximateIndustrial2021}.
\cite{Jungel2023} used a \gls*{CO}-enriched ML pipeline to learn a dispatching policy in a dynamic autonomous mobility-on-demand system.
While the work of~\cite{Jungel2023} may seem similar to this work at first sight, it differs in several fundamental aspects: the dispatching problem studied in \cite{Jungel2023} is solvable in polynominal time and hence does not rely on a heuristics in the \gls*{CO}-layer. Moreover the statistical model of \cite{Jungel2023} takes the form of a generalized linear model whereas we rely on deep learning in this work.

While \gls*{CO}-enriched ML pipelines are widely used to solve real-world problems, existing work in this field relies on exact algorithms in the \gls*{CO}-layer so far. However, this approach is not tractable in settings where the \gls*{CO} problem is difficult to solve, as is the case in most dynamic~\glspl*{VRP}. Moreover, existing work assumes CO-layers with linear objective functions where the predicted objective costs have the same dimension as the decision variables, which is not the case in our dynamic~\gls*{VRP} setting.

\subsection{Contributions}
To close the research gap outlined above, we propose a novel \gls*{ML}-based pipeline enriched with a \gls*{CO}-layer, and apply it to a novel class of dynamic~\glspl*{VRP} introduced in the EURO Meets NeurIPS Vehicle Routing Competition that is highly interesting for academia and practice. Specifically, our work contains several contributions. 
From a methodological perspective, we generalize the \gls*{CO}-layer of \gls*{CO}-enriched \gls*{ML} pipelines to non-linear objective functions.
Note that in this context, we also extend the open source library \texttt{InferOpt.jl} to such non-linear settings. Moreover, we present the first \gls*{CO}-enriched \gls*{ML} pipeline that utilizes a metaheuristic component to solve the \gls*{CO}-layer. By so doing, we show that our general \gls*{ML}-\gls*{CO} paradigm, which formally requires optimal solutions to derive true gradients, can work well with heuristic solutions in practice. In this context, we detail how to carefully design a metaheuristic that allows to derive heuristic solutions, i.e., approximate gradients, which enable convergence in the \gls*{ML}-layer of our pipeline. We show how to train the \gls*{ML}-layer of this pipeline in a supervised learning setting, i.e., based on a training set derived from an anticipative strategy.
We present a comprehensive numerical study to show the efficacy of our methodology in a benchmark against state-of-the-art approaches.
Beyond this, our study validates that our learning approach generalizes well to unseen scenarios. Interestingly, our results point at the fact that, counterintuitive to common practice, imitating an anticipative strategy can work well for high-dimensional multi-stage stochastic optimization problems. 

We refer to our git repository (\url{https://github.com/tumBAIS/euro-meets-neurips-2022}) for instructions and all material necessary to reproduce the results outlined in this paper.

\subsection{Outline}
The remainder of this paper is organized as follows. Section~\ref{sec:problem_setting} provides a formal definition of our problem setting before Section~\ref{sec:pipeline} introduces our \gls*{CO}-enriched \gls*{ML} pipeline. We then detail the individual layers of our pipeline. Specifically, Section~\ref{sec:combinatorial_problem} details the algorithmic framework used in the \gls*{CO}-layer, while Section~\ref{sec:learning} details the learning methodology for the \gls*{ML}-layer.
Section \ref{sec:results} details the design of our computational study and it's results. Finally, Section \ref{sec:conclusion} concludes the paper.

\section{Problem setting}\label{sec:problem_setting}

Our problem setting focuses on a variant of the dynamic \gls*{VRPTW} introduced in the EURO Meets NeurIPS Vehicle Routing Competition~\parencite[see][]{euromeetsneurips2022}.
In this dynamic \gls*{VRPTW}, we aim to find a cost-minimal set of routes that start and end at a central depot $\Depot$ and allow a fleet of vehicles to serve a set of requests $\FullSetOfRequests$ within a finite planning horizon $\PlanningHorizon$. We focus on an online problem setting in which the request set $\FullSetOfRequests$ is initially unknown and requests continuously arrive over $\PlanningHorizon$.
Within this planning horizon, the fleet operator makes dispatching decisions at (equidistant) time steps $\TimeStep \in \PlanningHorizon$ and needs to serve all requests that arrive during the planning horizon. Accordingly, we discretize the planning horizon into a set of $n$ epochs $\TimeHorizon =  \{[\TimeStep_0, \TimeStep_1], [\TimeStep_1, \TimeStep_2], \dots, [\TimeStep_{n-1}, \TimeStep_n]\}$ and denote the start time of an epoch $\Epoch$ as $\TimeStep_\Epoch$.

In each epoch, the fleet operator solves a dispatching and vehicle routing problem for the epoch dependent request set $\calR^e$.
Each request $\Request \in \calR^e$ has a certain demand $\Demand_\Request$, and vehicles have a homogeneous vehicle capacity $\VehicleCapacity$, which limits the maximum number of requests serviceable on a single route.
Serving a request $\Request$ takes $\ServiceTime_\Request$ time units, and a request must be served within a request-specific time window $\TimeWindow[\Request]$. Traveling from the delivery location of request $i$ to the delivery location of another request $j$ takes $\TravelTime_{ij}$ time units and incurs cost $\TravelCost_{ij}$. We can straightforwardly encode an instance of an epoch's planning problem 
on a fully-connected digraph $\EpochGraph = (\EpochVertices,\EpochArcs)$ with a vertex set $\EpochVertices = \calR^e \cup \{d\}$ where $d$ is the depot, and an arc set $\EpochArcs$.
With this notation, our problem representation unfolds as follows.

\paragraph{System state}
We describe a \emph{system state} at decision time $\TimeStep_\Epoch$ as $x^e = \EpochRequests$, where $e$ is the epoch starting at time $\TimeStep_\Epoch$.
Here, $\EpochRequests$ contains all requests revealed but not yet dispatched. 
In our specific setting no further information is required to describe the system state as there is no fleet limit.
Note that we use redundant notation $x^e = \EpochRequests$ to adhere to conventions commonly used in the domains of \gls*{ML} and vehicle routing, respectively.
We can distinguish requests contained in $\EpochRequests$ into two disjoint categories: \emph{must-dispatch} requests need to be dispatched in $e$ as dispatching these in a later epoch would violate the requests' time window upon delivery; \emph{postponable} requests can but do not have to be dispatched in $e$.

\paragraph{Feasible decisions}
Given the current state of the system $x^e$, the fleet operator chooses a subset of $\EpochRequests$ that will be served by vehicles leaving the depot in this epoch, and computes the respective routes that allow for vehicle dispatching. Vehicles dispatched in epoch $\Epoch$ leave the depot at time $\TimeStep_{\Epoch} + \DispatchTime$ and their routes cannot be modified once they have been dispatched, i.e., vehicles dispatched in epoch $\Epoch$ cannot serve requests revealed in epoch $\Epoch'$, with $\TimeStep_{\Epoch'} > \TimeStep_{\Epoch}$. In this context, a feasible decision $\decision^\Epoch \in \decisions(x^e)$ in state $x^e$ corresponds to a set of routes that 
\begin{itemize}
    \item[(i)] contains all \textit{must-dispatch} requests,
    \item[(ii)] allows each route to visit all contained requests within their respective time windows, and
    \item[(iii)] the cumulative customer demand on each route does not exceed the vehicle capacity~$\VehicleCapacity$. 
\end{itemize}
We can encode a feasible decision $\decision^\Epoch \in \calY(x^e)$ with a vector $(y^\Epoch_{i,j})_{(i,j) \in \EpochArcs}$ where
$$ y^\Epoch_{i,j} = \begin{cases}
1 & \text{if } (i,j) \text{ is in a route of the solution} \\
0 &  \text{otherwise.}
\end{cases}$$

\paragraph{System evolution}
The system transitions into the next epoch $\Epoch'$ once the fleet operator decides on $\decision^\Epoch \in \decisions(x^e)$. To describe $x^{\Epoch'}$, we derive $\mathcal{R}^{\Epoch'}$ by removing all requests contained in $y^\Epoch$ from $\EpochRequests$, and adding all requests that enter the system between $\TimeStep_{\Epoch}$ and $\TimeStep_{\Epoch'}$.

\paragraph{Policy}
Let~$\mathcal{X}$ denote the set of potential system states. Then, a (deterministic) policy~$\pi: \mathcal{X} \rightarrow \mathcal{Y}$ is a mapping that assigns a decision~$y^e \in \mathcal{Y}(x^e)$ to any system state~$x^e \in \mathcal{X}$.

\paragraph{Objective}
We aim to find a policy that minimizes the expected cost of serving all requests $\FullSetOfRequests$ over the planning horizon $\PlanningHorizon$. Formally,
\begin{equation}\label{eq:dynamicProblem}
    \min_{\policy} \mathbb{E} \left[ \sum_{e \in \SetOfEpochs} c(\policy(x^e)) \right],
\end{equation}
where $c: \decisions \rightarrow \mathbb{R}$ gives the cost of routes $\decision \in \decisions$.

\paragraph{Discussion}
The problem setting defined and formalized above contains various assumptions that might be questioned from a practitioner's perspective. In particular, assuming an unlimited fleet size and full knowledge of the request distribution appears to be rather unrealistic. As this paper focuses on the methodology used to win the EURO Meets NeurIPS Vehicle Routing Competition, we decided to keep these assumptions without further questioning for the sake of consistency and reproducibility. However, we like to emphasize that the methodological pipeline presented in this paper is readily applicable to problem settings with limited fleet sizes and incomplete knowledge of the underlying request distribution as long as some historical data is available.

\section{ML pipeline with CO-layer}\label{sec:pipeline}

\begin{figure}[t]
    \centering
    \includegraphics[width=\textwidth]{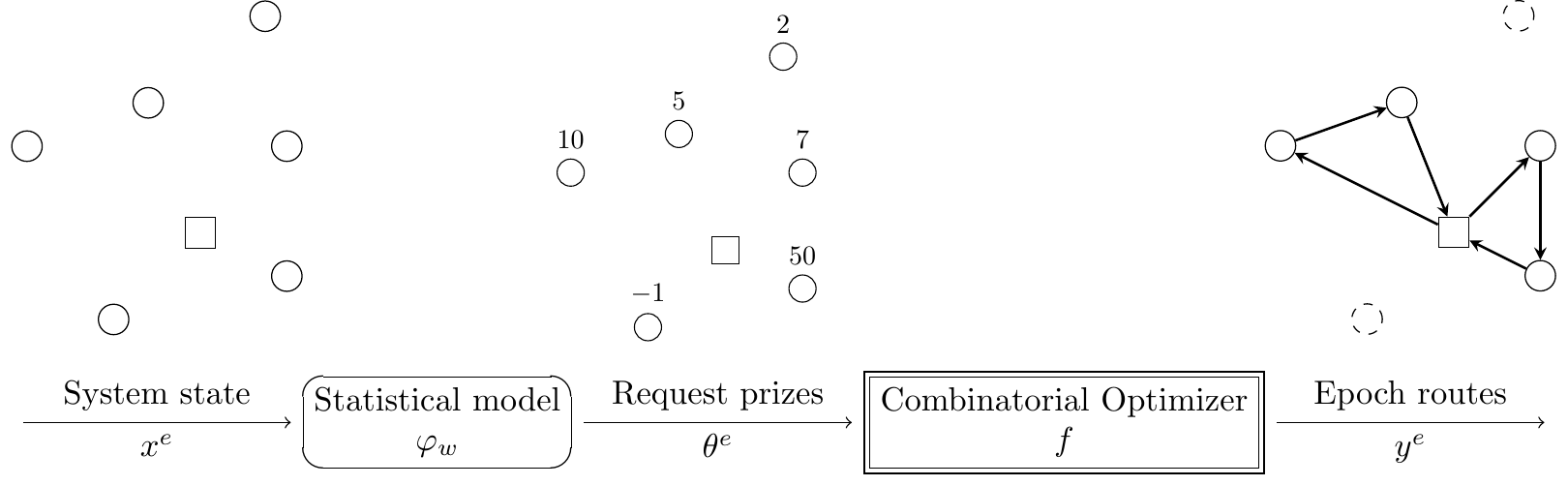}
    \caption{Our \gls*{CO}-enriched \gls*{ML} pipeline.}
    \label{fig:pipeline}
\end{figure}

To explain the rationale of our \gls*{CO}-enriched \gls*{ML} pipeline, we recall that a policy $\Policy$ maps a system state $x^e$ to a feasible decision $y^e$ in $\calY(x^e)$.
The state $x^e$ is a set of requests $\calR^e$, and the solution $y^e$ is a set of feasible decisions which encode the routes covering a subset of $\calR^e$.
The set of feasible decisions $\calY(x^e)$ hence coincides with the set of feasible solutions of a \acrshort*{PC-VRPTW}, a variant of the (static) \gls*{VRPTW} where it is not mandatory to serve all requests, but a prize $\theta^e_j$ is collected if request $j$ is served.
In this section, we show that any decision $y^e$ derived from an optimal policy in a dynamic problem setting corresponds to an optimal solution of a \acrshort*{PC-VRPTW} for a well chosen prize vector $\theta^e = (\theta^e_j)_{j \in \calR^e}$.
However, finding prizes $\theta^e_j$ is non-trivial, such that we resort to~\gls*{ML} for this purpose.
Specifically, we introduce a family of policies $(\pi_w)_{w}$ encoded by the~\gls*{CO}-enriched \gls*{ML} pipeline illustrated in Figure~\ref{fig:pipeline}:
in the \gls*{ML}-layer, a statistical model~$\varphi_w$ predicts $\theta^e$ based on the given system state $x^e$. This yields a \acrshort*{PC-VRPTW} instance $(x^e,\theta^e)$ which we solve in the \gls*{CO}-layer with a dedicated algorithm $f$. The algorithm's output $y^e$ then corresponds to our dispatching and routing decision.

In what follows, we first formally introduce the \acrshort*{PC-VRPTW} and proof that we can represent every optimal decision in epoch $e$ as a solution of a specially constructed \acrshort*{PC-VRPTW} instance. We then detail how we design our pipeline to leverage this observation for finding dispatching and routing decisions.

\paragraph{Prize-collecting VRPTW}
Recall that $\EpochGraph = (\EpochVertices,\EpochArcs)$ is a fully-connected digraph with vertex set $\EpochVertices = \calR^e \cup \{d\}$, comprising epoch requests $\calR^e$ and the depot $d$, and that a feasible solution of the \acrshort*{PC-VRPTW} $y^e \in \calY(x^e)$ can be encoded by the vector $(y_{i,j}^e)_{(i,j) \in \EpochArcs}$ where
$$ y_{i,j}^e = \begin{cases}
1 & \text{if } (i,j) \text{ is in a route of the solution} \\
0 &  \text{otherwise.}
\end{cases}$$
We further consider costs $(c_{i,j})_{(i,j) \in \EpochArcs}$ on each arc, and prizes $(\theta^e_j)_{j \in \calR^e}$ on each vertex.
Then, we can state the objective of the \acrshort*{PC-VRPTW} as follows,
\begin{equation}\label{eq:pc_vrptw}
    \max_{y \in \calY(x^e)} \underbrace{\sum_{\substack{(i,j) \in \EpochArcs\\ j\neq d} }\theta^e_j y_{i,j}}_{\text{total profit}} - \underbrace{\sum_{\substack{(i,j) \in \EpochArcs \\ \textcolor{white}{j} } } c_{i,j}y_{i,j}}_{\text{total routing cost}}.
\end{equation}

\begin{proposition}\label{prop:OptimalPrices}
For any $x^e$, there exists a $\theta \in \mathbb{R}^{|\calR^e|}$ such that any optimal solution of \eqref{eq:pc_vrptw} is an optimal decision with respect to~\eqref{eq:dynamicProblem}.
\end{proposition}
\proof{}
Since the horizon is finite and the set of feasible decisions at each step is also finite, there exists an optimal decision $y^\star$ for $x^e$.
Let $\bar \calR^e$ be the subset of requests of $\calR^e$ that are dispatched in $y^\star$. Then any solution $y$ which has lower or equally low routing costs and covers $\bar \calR^e$ exactly is also optimal.
This follows from the Bellman equation since the routes have no impact on the evolution of the state.
We can construct $y$ by solving a \acrshort*{PC-VRPTW} on $\calR^e$ with request prizes 
\begin{equation}
    \bar \theta_j = \begin{cases}
        M & \text{if } j \in \bar \calR^e\\
        -M & \text{otherwise,}
    \end{cases}
\end{equation}
where $M = \left(|\calR^e| \cdot \max\limits_{(i,j)\in \EpochArcs}c_{i,j}\right)$ is a large constant.
The corresponding \acrshort*{PC-VRPTW} solution $y$ clearly covers $\bar \calR^e$ exactly and has at most the routing cost of $y^\star$.
\hfill $\square$
\endproof

\paragraph{\Gls*{CO}-layer}
We embed the \acrshort*{PC-VRPTW} as a \emph{layer} in our \gls*{CO}-enriched \gls*{ML} pipeline. Hence, this \gls*{CO}-layer must support forward and backward passes to assure compatibility with the \gls*{ML}-layer.

The forward pass simply solves the \acrshort*{PC-VRPTW} instance defined by $(x^e, \EpochPrizes)$ using a metaheuristic algorithm $f$ detailed in Section~\ref{sec:combinatorial_problem}.
The backward pass backpropagates the gradient of the loss used in the learning algorithm through this layer.
Section~\ref{sec:learning} introduces this loss and its gradient. To make their statement easier, we reformulate~\eqref{eq:pc_vrptw} as
\begin{equation}\label{eq:co_layer}
    f\colon \theta^e \mapsto\argmax_{y \in \calY(x^e)} \; {\theta^e}^{\top} g(y) + h(y) \quad \text{where} \quad g(y) = \bigg(\sum_{i\in \mathcal{V}^e} y_{i,j}\bigg)_{j\in \calR^e} \text{ and } h(y) = \sum_{(i,j) \in \EpochArcs}c_{i,j}y_{i,j}.
\end{equation}

\paragraph{ML-layer}
\label{sec:predictor}
Finding the optimal prizes $\EpochPrizes$ of Proposition~\ref{prop:OptimalPrices} is non-trivial.
We therefore use a statistical model~$\varphi_w$ to predict a vector of prizes $\EpochPrizes = \predictor\left(\EpochInstance\right) \in \mathbb{R}^{|\EpochRequests|}$ given the system state~$x^e$.
The only technical aspect from the \gls*{ML} perspective is that the dimensions of the input and the output are not fixed. Indeed, the number of requests may change from one state to another, and also differs across instances.
We benchmark different choices for commonly used statistical models~$\varphi_w$ in our computational study (cf. Section~\ref{sec:results}). 

\paragraph{Discussion}
Alternative \gls*{ML}-based solution approaches presented in the EURO Meets NeurIPS Vehicle Routing Competition \parencite[e.g.,][]{optiml2022} generally proceed in two steps.
They first apply a binary classifier, to decide on which requests to dispatch and postpone, respectively. They then construct routes covering these dispatched requests in a second step, essentially decoupling request dispatching from route construction.
A major difficulty in this approach is to learn a statistical model which implicitly balances the current route costs and future route costs to find optimal dispatching decisions.
We bypass this difficulty by taking dispatching and routing decisions simultaneously in our \gls*{CO}-layer~\eqref{eq:co_layer}, which allows us to train our statistical model based on the routing rather than the dispatching decision.

\section{Combinatorial optimization algorithm}
\label{sec:combinatorial_problem}

To derive solutions of the~\acrshort*{PC-VRPTW} within our \gls*{CO}-layer, we propose a metaheuristic algorithm based on \gls*{HGS}.
Our algorithm extends the implementation of the \gls*{HGS} algorithm introduced by \cite{koolHybrid2022}, which adapts the original \gls*{HGS} of \cite{vidalHybrid2021} to the \gls*{VRPTW}.

Tailoring the work of \cite{koolHybrid2022} to our problem setting requires various modifications to support optional requests, e.g., adaptions of the local search, initialization, and crossover procedures. We further introduce new mutation mechanisms, implement \acrshort*{PC-VRPTW} specific neighborhoods, and allow to warm-start the population. In what follows, we summarize the core concepts of the \gls*{HGS} algorithm before detailing our \acrshort*{PC-VRPTW}-specific modifications.

\begin{figure}[h]
    \centering
    \includegraphics[width=\textwidth]{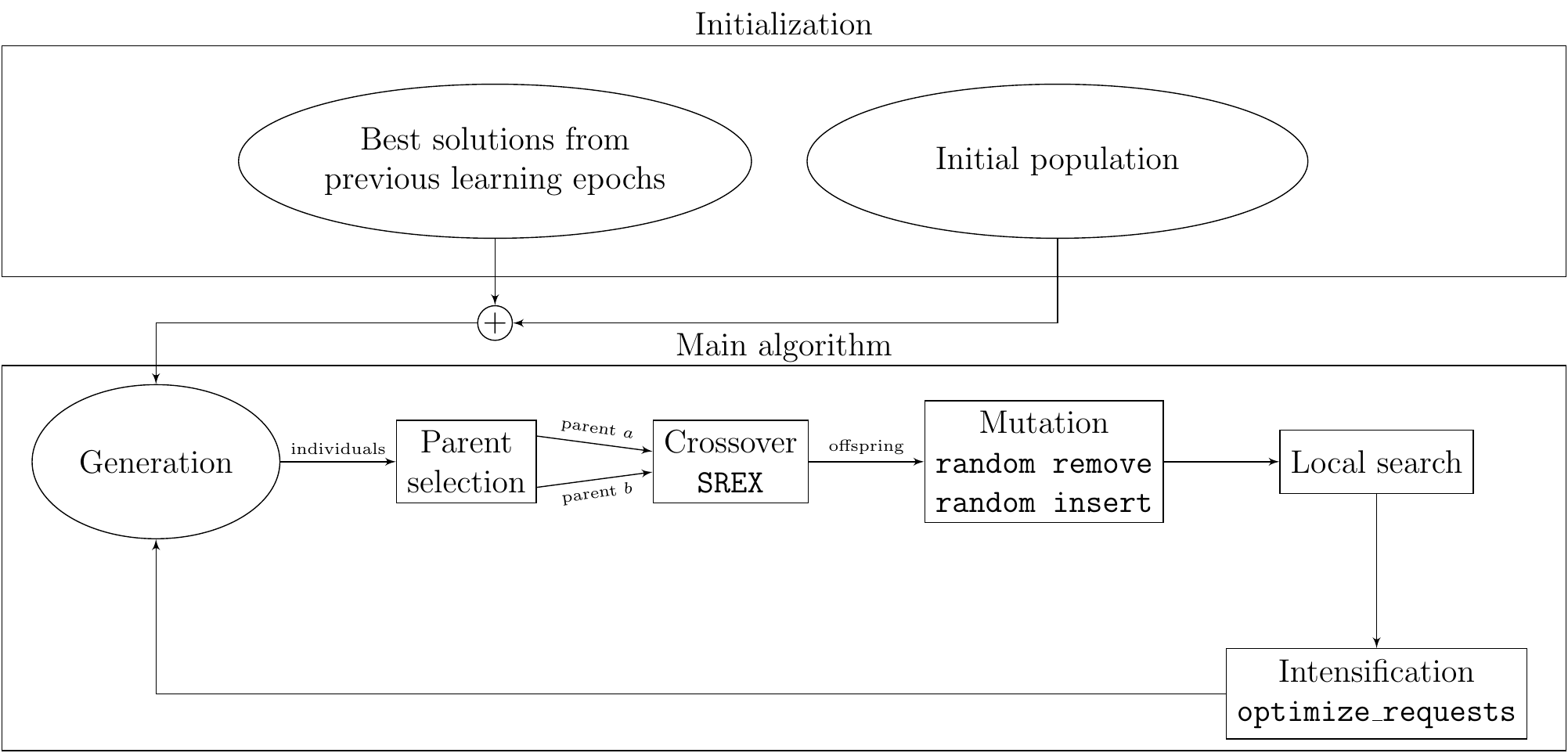}
    \caption{General structure of our metaheuristic algorithm. Round nodes indicate populations, square nodes indicate algorithmic components.}
    \label{figure:pchgs}
\end{figure}

\paragraph{Hybrid Genetic Search for the Vehicle Routing Problem With Time Windows}
The \gls*{HGS} algorithm is an evolutionary algorithm that maintains a \textit{population} of solutions, organized in two disjoint sub-populations that contain feasible and infeasible solutions, respectively.
The algorithm makes this population evolve over time, generating \textit{offspring} solutions by combining promising \textit{parent} solutions selected from the population in a randomized binary tournament.
The algorithm then improves the generated offspring solution in a local search procedure. Note that it uses a penalty-based approach to explore infeasible regions of the solution space. The local search procedure yields a locally optimal solution, which is then added to the population. This may trigger a \textit{survivor selection} procedure if the population size exceeds a certain threshold. This procedure eliminates a subset of solutions from the population. Survivor and parent selection are based on the \textit{fitness} of a solution, a metric which captures the quality, i.e., objective value, of a solution and it's contribution to the population's diversity. This ensures a sufficiently diverse population, balancing diversification and intensification in the genetic algorithm. To utilize this algorithmic structure for our problem setting, i.e., the \acrshort*{PC-VRPTW}, we applied the following adaptions and extensions, leading to the algorithm outlined in Figure~\ref{figure:pchgs}.

\paragraph{Solution representation} We represent decisions on which requests to serve and which to ignore in our solution representation implicitly. Specifically, we use the same giant-tour representation as in \cite{vidalHybrid2021}, but allow incomplete giant-tours. Here, requests deemed unprofitable are absent from the giant-tour and thus not considered by traditional local search operators.
\paragraph{Accounting for optional requests during crossover}
The \gls*{HGS} proposed in \cite{koolHybrid2022} generates offspring solutions using two crossover operators: \gls*{OX} \parencite{oxCrossover} and \gls*{SREX} \parencite{nagataSrex}. As a first modification, we remove the \gls*{OX} operator as our benchmarking experiments indicate that it has no substantial impact on the algorithm's performance for our problem variant. Our second modification amends the \gls*{SREX} operator to the prize-collecting \gls*{VRPTW}. Specifically, we preserve the set of requests served by the first parent in any generated offspring by re-inserting requests that are currently not served as in regular \gls*{SREX}, i.e., sequentially. 

\paragraph{Additional diversification and intensification mechanisms}
We further introduce two new mutation operators that diversify and intensify solutions based on the set of served requests.
The first operator (\texttt{random remove/insert}), based on the result of a coin toss, either removes $\ParamMutationRemovalFactor \cdot |\RequestSet|$ of the requests currently served, or inserts $\ParamMutationInsertionFactor \cdot |C \setminus \RequestSet|$ of the currently unserved requests. This mutation occurs with a probability of $\ParamMutationProbability$ right after offspring generation.

The second operator (\texttt{optimize\_request\_set}) optimizes the set of requests served by a given solution.
Specifically, this operator first removes any requests from the solution that cause a detour whose cost is higher than the request's profit, and re-inserts any profitable requests that are not part of the current solution.
The operator perturbs insertion and removal costs with a random factor drawn uniformly from the interval $[\ParamOptimizeRequestSetPerturbationLB, \ParamOptimizeRequestSetPerturbationUB]$. In contrast to the first operator, we delay evaluating the second operator until the LS converges to avoid removing an excessive amount of requests due to poor solution quality. We run this operator with probability $\ParamOptimizeRequestSetProbability$. To intensify our search in regions around promising solutions, we further run this operator (without perturbation) on any solution that improves on the current best solution.

\paragraph{Local Search}
Table~\ref{table:operators} shows the local search operators used in our algorithm. We consider traditional and prize-collecting \gls*{VRPTW} specific operators. Specifically, traditional operators refer to those defined for the \gls*{VRPTW} \parencite[see][]{koolHybrid2022, vidalHybrid2021}. These traditionally explore neighborhoods by exchanging arcs within a solution, i.e., by changing the position of one or several requests in the current solution. Hence, they preserve the set of requests that receive service in a given solution. Prize-collecting \gls*{VRPTW} specific operators on the other hand work with the request set exclusively, that is, they remove or insert a request from a given solution. We evaluate these after evaluating traditional \gls*{VRPTW} operators to avoid removing profitable requests prematurely, i.e., due to bad routing decisions.
As in \cite{vidalHybrid2021}, we further distinguish between small and large neighborhood operators: small neighborhoods consider only those moves that involve requests which are geographically close and compatible w.r.t. their time windows. Large neighborhoods on the other hand remain unrestricted.

\begin{table}[h]
    \footnotesize
    \begin{threeparttable}
    \begin{tabularx}{\textwidth}{p{2cm} p{3cm} X}
    \hline
    Type & Name & Description \\ \hline
    Traditional (Small) & \texttt{relocate} & removes a single request from it's route and re-inserts it at a different position in the solution \\
     & \cellcolor{LightGray} \texttt{relocate pair}  & \cellcolor{LightGray} removes a pair of consecutive requests from their route and re-inserts them at a different position in the solution \\
     & \texttt{relocate reversed pair} & removes a pair of consecutive requests from their route and re-inserts them in reverse order at a different position in the solution \\
     & \cellcolor{LightGray} \texttt{swap} & \cellcolor{LightGray} exchanges the position of two requests in the solution \\
     & \texttt{swap pair} & exchanges the positions of two pairs of consecutive requests in the solution \\
     & \cellcolor{LightGray} \texttt{swap pair with single} & \cellcolor{LightGray} exchanges the positions of a pair of consecutive requests with a single request in the solution \\
     & \texttt{2-opt} & reverses a route segment \\
     & \cellcolor{LightGray} \texttt{2-opt*} & \cellcolor{LightGray} splits two routes into two segments each, swapping a segment from the first with a segment from the second route \\
    Traditional (Large) & \texttt{relocate*} & removes a request from it's route and inserts it at the best possible position in any spatially overlapping route \\
     & \cellcolor{LightGray} \texttt{swap*} & \cellcolor{LightGray} removes requests $a$ and $b$ from routes $r_a$ and $r_b$ with spatial overlap, inserting $a$ and $b$ at the best position in $r_b$ and $r_a$, respectively \\
    PC-VRPTW & \texttt{serve request} & inserts a currently unserved request into the solution \\
     & \cellcolor{LightGray} \texttt{remove request} & \cellcolor{LightGray} removes a request from the solution \\ \hline
    \end{tabularx}
    \end{threeparttable}
    \caption{Local search operators}
    \label{table:operators}
\end{table}

\paragraph{Initialization}
We apply a pre-processing technique to account for extreme weights encountered during the learning procedure. Specifically, we maintain a set of certainly profitable and certainly unprofitable requests based on the maximum and minimum detour required, respectively. We determine certainly unprofitable requests based on the following observation: if $\min_{j \in \EpochRequests \cup \{\Depot\}} \TravelCost_{jr} + \min_{j \in \EpochRequests\cup \{\Depot\}} \TravelCost_{rj} - \Prize_r \geq \max_{i, j \in (\EpochRequests\cup \{\Depot\})^2} c_{i, j}$ holds for some request $r \in \EpochRequests$, then for any route $\sigma$ that includes $r$, a cheaper route $\sigma'$ exists. To see this, let $u, v \in \EpochRequests \cup \{\Depot\}$ such that $(u, r), (r, v) \in \sigma$. Then $c_{u,r} + c_{r, v} - \Prize_r \geq \max_{i, j \in (\EpochRequests\cup \{\Depot\})^2} c_{i, j} \geq c_{u, v}$, such that removing $r$ from $\sigma$ yields a cheaper route. Analogously, we force the inclusion of a request $r \in \EpochRequests$ if $\TravelCost_{d,r} + \TravelCost_{r,d} < \Prize_r$ holds. This procedure further allows to account for must-dispatch customers during the evaluation phase by setting $\Prize_r$ accordingly.

The magnitude of changes in customer prices decreases as the learning procedure converges. We exploit this behaviour in our algorithm and seed the initial population, generated as in \cite{koolHybrid2022}, with solutions of previous learning epochs. Specifically, we generate promising solutions in a new construction heuristic which first applies our \texttt{optimize\_customers} operator to the solution to account for changed request prizes, and then optimizes the routing decisions using the local search procedure. Note that this is only possible while training the pipeline. When evaluating the pipeline, we use starting solutions generated as in \cite{koolHybrid2022} exclusively.

\paragraph{Discussion}
Focusing on our algorithmic design decisions, we note that we have tailored our algorithm to the unique challenges of our learning methodology and training environment. Specifically, successfully computing approximate gradients that ensure convergence during training requires solutions with low variance. 
Here, the algorithm needs to generalize well to the different prize distributions observed during training. Beyond this, it must be capable of providing such solutions within tight time constraints to allow training within a reasonable amount of time.

We have tackled these challenges with a design that focuses on guiding the population towards promising regions of the search space through aggressive diversification (e.g., randomized insertion and deletion of requests) and intensification (\texttt{optimize\_request\_set}) operators, explicit request insertion and deletion neighborhoods, and warmstarting. 
The resulting algorithm behaves more greedily than implicit approaches \parencite[e.g.][]{vidalPrizeCollecting}, but converges reliably within the tight time constraints imposed during training. 
We further cope with extreme prizes observed during training in a preprocessing step.

\section{Learning approach}
\label{sec:learning}
The objective of our learning problem is to find parameters~$w$ such that the statistical model $\Predictor$ predicts a prize vector~$\theta$ that leads to ``good'' decisions in the \gls*{CO}-layer.
To reach this objective, we train our \gls*{CO}-enriched \gls*{ML} pipeline to imitate a good policy.
To do so, we follow a supervised learning setting
and therefore build a dataset $\mathcal{D} = \{(\EpochInstance[1], \TargetEpochRoutes[1]), \dots, (\EpochInstance[n], \TargetEpochRoutes[n])\}$ of state instances $\EpochInstance[i]$ with the decisions $\TargetEpochRoutes[i]$ taken by the imitated policy. In this supervised learning setting, we define a loss $\mathcal{L}(\theta, \bar y)$ which quantifies the error when we predict $f(\theta)$ instead of $\bar y$. Then, we formulate the learning problem as finding the parameter $\hat w_\mathcal{D}$ that minimizes the empirical risk, i.e., the average loss on the training dataset,
\begin{equation}
    \label{eq:learning-problem}
    \hat\Weights_{\mathcal{D}} = \argmin_\Weights \sum_{i=1}^n \mathcal{L}(\varphi_\Weights(\EpochInstance[i]), \TargetEpochRoutes[i]).
\end{equation}

The rest of this section introduces the imitated policy, the training set, the loss, and our algorithm to solve learning problem~\eqref{eq:learning-problem}. 

\subsection{Anticipative decisions in the training set}
\label{subsec:anticipative_policy}

The dynamic~\gls*{VRPTW} is difficult because future requests are unknown when taking decisions.
If we know the future, i.e., if we knew all the requests at the beginning of the horizon, we would just have to solve the corresponding static VRPTW at the beginning of the horizon, and then take the corresponding dispatching decisions at each epoch.
This procedure would give us the optimal \text{anticipative} policy.
However, such an anticipative policy can of course not be used in practice because it relies on unavailable information.
Nonetheless, we can rebuild the decisions taken by the anticipative policy a posteriori, when all information has been revealed, i.e., on historical data.
We therefore can learn to imitate an anticipative policy, which we gained from historical data, in our learning problem~\eqref{eq:learning-problem}.

Practically, we rebuild the decision of the anticipative policy as follows.
We know all the requests of an historical \gls*{VRPTW} instance.
Then, we associate a \textit{release time} to each request $i$, which is equal to the time $\tau_i$ at which the request would be revealed in a dynamic setting, so the earliest point in time at which request $i$ can be served.
We then seek routes to serve all these requests at minimum cost while respecting the release times and the time windows. We use the adapted \gls*{HGS} of~\cite{koolHybrid2022} to solve this variant of the static \gls*{VRPTW}.
This yields a set of routes $P$ from which we reconstruct the decisions taken in each epoch as follows. 
For route $p \in P$, let $\TimeStep_p = \max_{r \in p} \TimeStep_{r}$ be the latest release time of any request route $p$ serves. The anticipative policy dispatches route $p$ in the first epoch $e$ where $\TimeStep_e \geq \TimeStep_p$, i.e., the first epoch where all served requests have been revealed.

\subsection{Loss function}\label{subsec:loss_function}
We recall that the objective of our learning problem is to find parameter values $w$ of a statistical model~$\varphi_w$ such that for any state $x$, the CO-layer predicts a good decision~$y=f(\theta)$ for the \acrlong*{PC-VRPTW} instance~$(x, \theta=\varphi_w(x))$.
More precisely, for each state-decision pair $(x, \bar y)$ in the data set, we want the target decision~$\TargetEpochRoutes[]$ to be as close as possible to the optimal solution of the \acrlong*{PC-VRPTW} \eqref{eq:pc_vrptw}.
Hence, it is natural to take the non-optimality of $\bar y$ as a solution of~$f(\theta)$ as loss function:
\begin{equation}
    \label{eq:natural-loss}
    \Loss(\Prize, \TargetEpochRoutes[]) = \max_{\EpochRoutes[]\in\decisions(\EpochInstance[])} \{\Prize^\top g(\EpochRoutes[]) + h(\EpochRoutes[])\} - (\Prize^\top g(\TargetEpochRoutes[]) + h(\TargetEpochRoutes[])).
\end{equation}
We clearly have $\calL(\theta,\bar y) \geq 0$ in general, and $\calL(\theta,\bar y) = 0$ only if $y$ is an optimal solution of~\eqref{eq:co_layer}.
Unfortunately, the polyhedron $\mathcal{P}(\bar y) = \{\theta\in\mathbb{R}^{|x|} \,|\, \theta^\top g(\bar y) + h(\bar y) \geq \theta^\top g(y) + h(y),\, \forall y\in\mathcal{Y}(x)\}$ is generally highly degenerate.
For instance, if $h = 0$, then $\theta = 0$ belongs to $\calP(y)$ for all $y$ in $\calY$.
This implies that $\theta = 0$ would be an optimum of our learning problem, which is a problem because such a $\theta$ allows our \gls*{CO}-layer to return any solution in $\calY(x)$.
This degeneracy can be removed by considering the perturbed loss
\begin{equation}
    \label{eq:perturbed-loss}
    \mathcal{L}_\varepsilon(\theta,\bar y) = \mathbb{E}\Big[\max_{\EpochRoutes[]\in\decisions(\EpochInstance[])} \{(\Prize + \varepsilon Z)^\top g(\EpochRoutes[]) + h(\EpochRoutes[])\}\Big] - (\Prize^\top g(\TargetEpochRoutes[]) + h(\TargetEpochRoutes[])).
\end{equation}
where $\varepsilon > 0$ and $Z \sim \mathcal{N}(0,\mathbf{I}_{|x|})$ is a Gaussian perturbation.
The perturbed loss $\mathcal{L}_\varepsilon$ has been considered in the literature only in the case where $g(y) = y$ and $h(y)=0$~\parencite[cf.][]{berthetLearning2020}. 
Due to non-zero $h$, the geometry of the loss changes. 
Notably, while the size of $\varepsilon$ does not matter when $h = 0$ \parencite[cf.][]{parmentierLearning2021-1}, it does matter when $h$ is non-zero: the larger $\varepsilon$, the smaller the impact of $h$ on the prediction.
Proposition~\ref{prop:lossProperties} summarizes the geometry of these losses in our more general context.
\begin{proposition}\label{prop:lossProperties}
Let $x\in\mathcal{X}$, $\bar y\in\mathcal{Y}(x)$.
Let $\mathcal{C}(\bar y) = \{\theta\in\mathbb{R}^{|x|} \colon \theta^\top g(\bar y) \geq \theta^\top g(y),\,\forall y\in\mathcal{Y}(x)\}$ be the \emph{normal cone} associated to $g(\bar y)$.

\begin{enumerate}
    \item $\theta \mapsto \calL(\theta,\bar y)$ is piecewise linear and convex, with subgradient $$\underbrace{g\Big(\argmax_{\EpochRoutes[]\in\decisions(\EpochInstance[])} \theta ^\top g(\EpochRoutes[]) + h(\EpochRoutes[])\Big)}_{=g(f(\theta))} - g(\bar y) \in \partial_\theta \calL(\theta,\bar y).$$
    \item $\theta \mapsto \mathcal{L}_\varepsilon(\theta,\bar y)$ is $\calC^{\infty}$ and convex with gradient $$\nabla_\theta  \mathcal{L}_\varepsilon(\theta, \bar y) = \mathbb{E} \left[g\left(\argmax_{\EpochRoutes[]\in\decisions(\EpochInstance[])} (\theta + \varepsilon Z)^\top g(\EpochRoutes[]) + h(\EpochRoutes[])\right)\right] - g(\TargetEpochRoutes[]) = \mathbb{E}[g(f(\theta+\varepsilon Z))] - g(\bar y).$$
    \item $\mathcal{L}_\varepsilon(\theta,\bar y) \geq \mathcal{L}(\theta,\bar y)$.
    \item $\mathcal{C}(\bar y)$ is the recession cone of $\calP(\bar y)$.
    \item Let $\theta\in\mathbb{R}^{|x|}$. If $\eta$ is in $\mathcal{C}(\bar y)\backslash\{0\}$, then $\lambda \mapsto \calL(\theta + \lambda \eta, \bar y)$ is non increasing. If in addition $\mathcal{C}(\bar y)\neq\mathbb{R}^{|x|}$, then $\lambda \mapsto \mathcal{L}_\varepsilon(\theta + \lambda \eta, \bar y)$ is decreasing. \label{prop:item-decreasing}
    \item Let $\theta\in\mathbb{R}^{|x|}$. If $\eta$ is in the interior $\mathring{\mathcal{C}}(\bar y)$ of $\mathcal{C}(\bar y)$, then $ \displaystyle\lim_{\lambda \rightarrow \infty}\mathcal{L}(\theta + \lambda \eta, \bar y) = \lim_{\lambda \rightarrow \infty}\mathcal{L}_\varepsilon(\theta + \lambda \eta, \bar y) = 0$.
\end{enumerate}
\end{proposition}

\proof{}
\begin{enumerate}
    \item $\theta \mapsto \calL(\theta,\bar y)$ is linear and convex as it is the maximum of mappings that are linear in~$\theta$. 
    Let $(\theta, \tilde\theta)\in\mathbb{R}^{|x|}\times\mathbb{R}^{|x|}$.
    Let $y^*$ be in $\argmax_{\EpochRoutes[]\in\decisions(\EpochInstance[])}\theta ^\top g(\EpochRoutes[]) + h(\EpochRoutes[])$ and $\tilde y$ in $\argmax_{\EpochRoutes[]\in\decisions(\EpochInstance[])}\tilde\theta ^\top g(\EpochRoutes[]) + h(\EpochRoutes[])$. By definition of $\tilde y$, we have $\tilde\theta ^\top g(\tilde y) + h(\tilde y) \geq \tilde\theta ^\top g(y^*) + h(y^* )$, which gives $ \calL(\tilde \theta,\bar y) \geq \calL(\theta^*,\bar y) + (\tilde \theta - \theta^*)^\top \big(g(y^*)-g(\bar y)\big)$ and the subgradient in $\partial_\theta\mathcal{L}(\theta, \bar y)$.
    \item Since $\mathcal{L}_\varepsilon(\theta,y) = \mathbb{E}[\calL(\theta + \varepsilon Z,y)]$, we obtain that $\theta \mapsto \mathcal{L}_\varepsilon(\theta,y)$ is $C^{\infty}$ as a convolution product of $\theta\mapsto\mathcal{L}(\theta, \bar y)$ with a Gaussian density (which is $C^{\infty}$), and is convex as an expectation of a convex function. The gradient follows from the subgradient of $\calL(\theta,y)$.
    \item Let $\theta\in\mathbb{R}^{|x|}$. For all $\tilde y$ in $\mathcal{Y}(x)$, we have
    \begin{align*}
        \text{for all } z\in\mathbb{R}^{|x|},\quad \max_y (\theta + \varepsilon z)^\top g(y) + h(y) &\geq (\theta + \varepsilon Z)^\top g(\tilde y) + h(\tilde y)\\
        \text{hence,}\quad \mathbb{E}[\max_y (\theta + \varepsilon Z)^\top g(y) + h(y)] &\geq \theta^\top g(\tilde y) + h(\tilde y)\quad(\text{since }\mathbb{E}[Z] = 0)\\
        \text{and}, \quad \mathbb{E}[\max_y (\theta + \varepsilon Z)^\top g(y) + h(y)] - (\theta^\top g(\bar y) + h(\bar y)) &\geq \theta^\top g(\tilde y) + h(\tilde y) - (\theta^\top g(\bar y) + h(\bar y)).
    \end{align*}
    Finally: $\mathbb{E}[\max_y (\theta + \varepsilon Z)^\top g(y) + h(y)] - (\theta^\top g(\bar y) + h(\bar y)) \geq \max_y \{\theta^\top g(y) + h(y)\} - (\theta^\top g(\bar y) + h(\bar y)),$
    i.e. $\boxed{\mathcal{L}_\varepsilon(\theta,\bar y) \geq \mathcal{L}(\theta,\bar y)}.$
    \item Let $\theta\in\mathcal{P}(\bar y)$. We have $\forall y,\, (\theta + \lambda \eta)^\top g(\bar y) + h(\bar y) \geq (\theta + \lambda \eta)^\top g(y) + h(y) \Leftrightarrow \forall y,\, \theta^Tg(\bar y) + h(\bar y) + \lambda \eta^\top \big(g(\bar y) - g(y)\big) \geq \theta^\top g(y) + h(y)$.
    If $\eta$ is in $\calC(\bar y)$, then $\eta^\top \big(g(\bar y) - g(y)\big) \geq 0$, and $\theta + \lambda \eta$ is in $\calP(\bar y)$ for any $\lambda > 0$. If $\eta$ is not in $\calC(\bar y)$, then $\eta^\top \big(g(\bar y) - g(y)\big) < 0$ and there exists a $\lambda > 0$ such that $\theta + \lambda \eta$ is not in $\calP(\bar y)$.
    \item Let $\eta\in\calC(\bar y)$. Let us denote by $\hat y = \argmax_y (\theta + \lambda\eta)^\top g(y) + h(y)$.
    We have
    $$\text{for all }\lambda,\quad \calL(\theta + \lambda\eta, \bar y) = (\theta + \lambda\eta)^\top g(\hat y) + h(\hat y) - [(\theta + \lambda\eta)^\top g(\bar y) + h(\bar y)]$$
    By taking the derivative with respect to $\lambda$, we obtain: $\eta^\top (g(\hat y) - g(\bar y))$ which is negative by definition of $\eta$. This gives us the non increasing-property of $\lambda\mapsto\calL(\theta + \lambda\eta, \bar y)$.

    Similarly, by denoting $\hat y_\varepsilon(Z) = \argmax_y (\theta + \lambda\eta + \varepsilon Z)^\top g(y) + h(y)$ for all $Z$, we obtain the derivative of $\lambda\mapsto\calL_\varepsilon(\theta + \lambda\eta)$: $\eta^\top(\mathbb{E}[g(\hat y_\varepsilon(Z))] - g(\bar y))$.
    If $\mathcal{C}(\bar y)\neq\mathbb{R}^{|x|}$, then $\mathbb{P}(\theta+\lambda\eta+\varepsilon Z \notin \mathcal{C}(\bar y)) > 0$. Hence, $\mathbb{E}[g(\hat y_\varepsilon(Z)]\neq g(\bar y)$, $\eta^\top(\mathbb{E}[g(\hat y_\varepsilon(Z)] - g(\bar y)) < 0$, and therefore $\lambda\mapsto\calL_\varepsilon(\theta + \lambda\eta)$ is decreasing.
    \item Let $\eta\in\mathring{\mathcal{C}}(\bar y)$. By definition of $\mathring{\mathcal{C}}(\bar y)$, we have $\eta^\top (\bar y - y) > 0$ for all $y\neq \bar y$.
    Hence, there exists $M>0$ such that for all $\lambda \geq M \text{ and } y\neq \bar y$, we have $(\theta + \lambda\eta)^\top\bar y + \geq (\theta + \lambda\eta)^\top y $.
    That is, for all $\lambda\geq M$ we have $\theta + \lambda\eta\in\calP(\bar y)$, i.e. $\mathcal{L}(\theta + \lambda\eta, \bar y) = 0$. Hence $\boxed{\lim\limits_{\lambda\rightarrow\infty}\mathcal{L}(\theta + \lambda\eta, \bar y) = 0}$.

    If $\mathcal{C}(\bar y)\neq\mathbb{R}^{|x|}$, from point~\ref{prop:item-decreasing} and the previous limit, we get that $\lambda \mapsto \mathcal{L}(\theta + \varepsilon z + \lambda\eta, \bar y)$ monotonically decreases to zero for any $z$.
    The monotone convergence theorem therefore gives $\boxed{\lim\limits_{\lambda\rightarrow\infty}\mathcal{L}_\varepsilon(\theta + \lambda\eta, \bar y) = 0}$. If $\mathcal{C}(\bar y)=\mathbb{R}^{|x|}$, the result is immediate.
\end{enumerate}
\endproof

\begin{figure}[H]
    \centering
    \begin{subfigure}[b]{0.46\textwidth}
        \centering
        \includegraphics[width=\textwidth, trim={90 60 12 40}, clip]{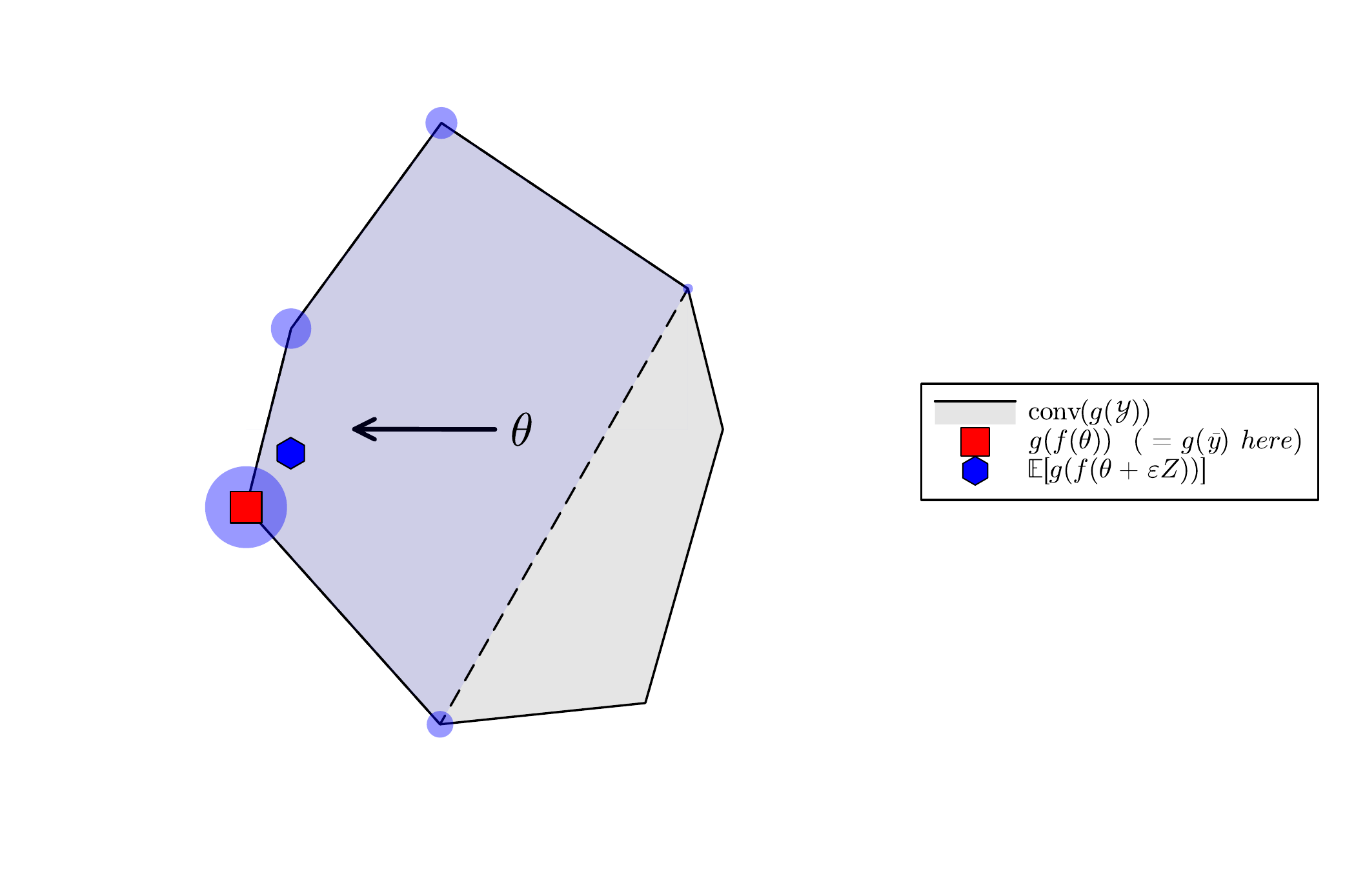}
        \caption{Polytope $\conv(g(\calY))$}
        \label{subfig:polytope}
    \end{subfigure}
    
    \begin{subfigure}[b]{0.45\textwidth}    
        \includegraphics[width=\textwidth]{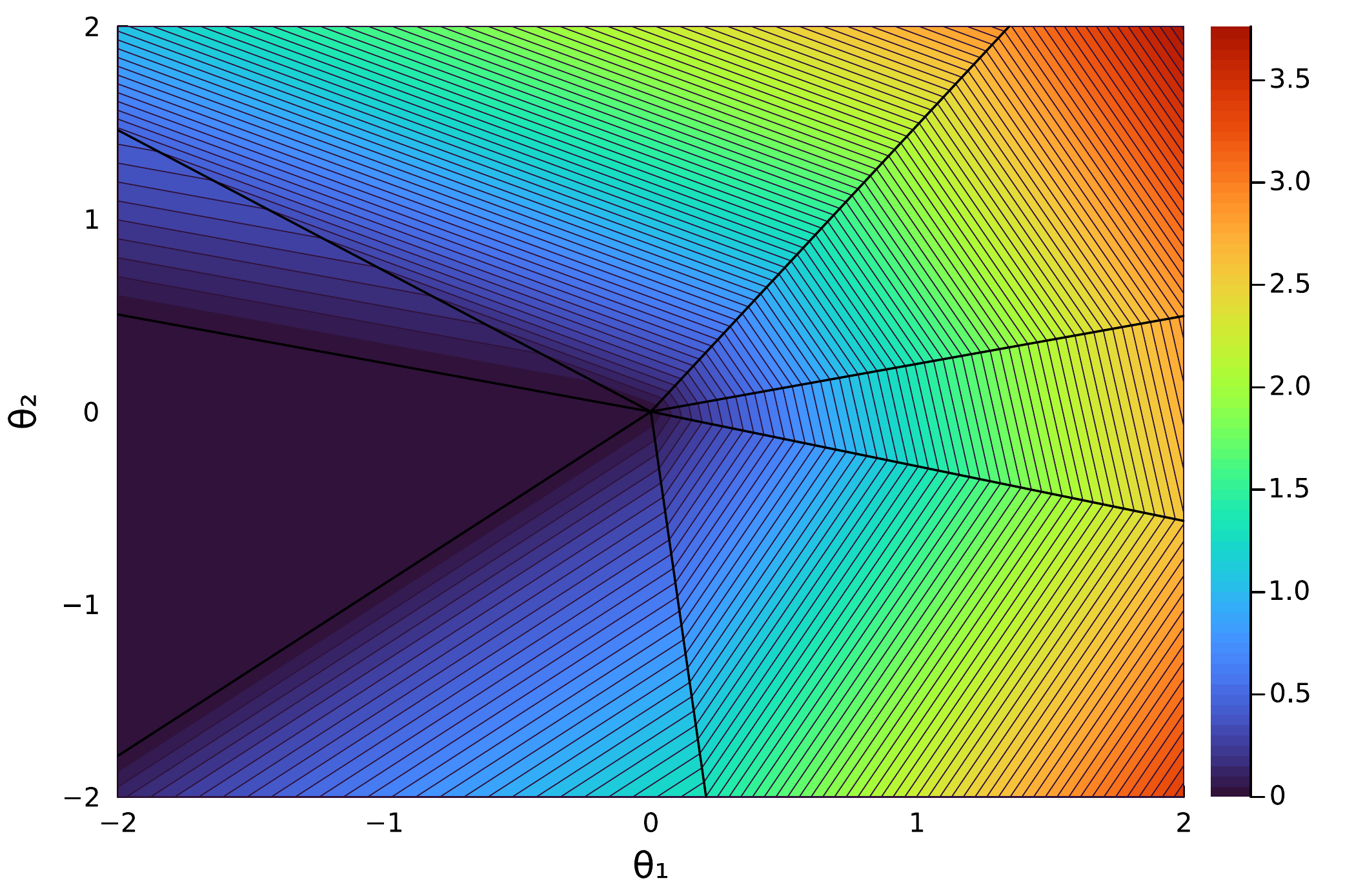}
        \caption{$\theta\mapsto\calL(\theta,\bar y)$ for $h=0$}
        \label{subfig:FYL}
    \end{subfigure}
    \begin{subfigure}[b]{0.45\textwidth}
        \includegraphics[width=\textwidth]{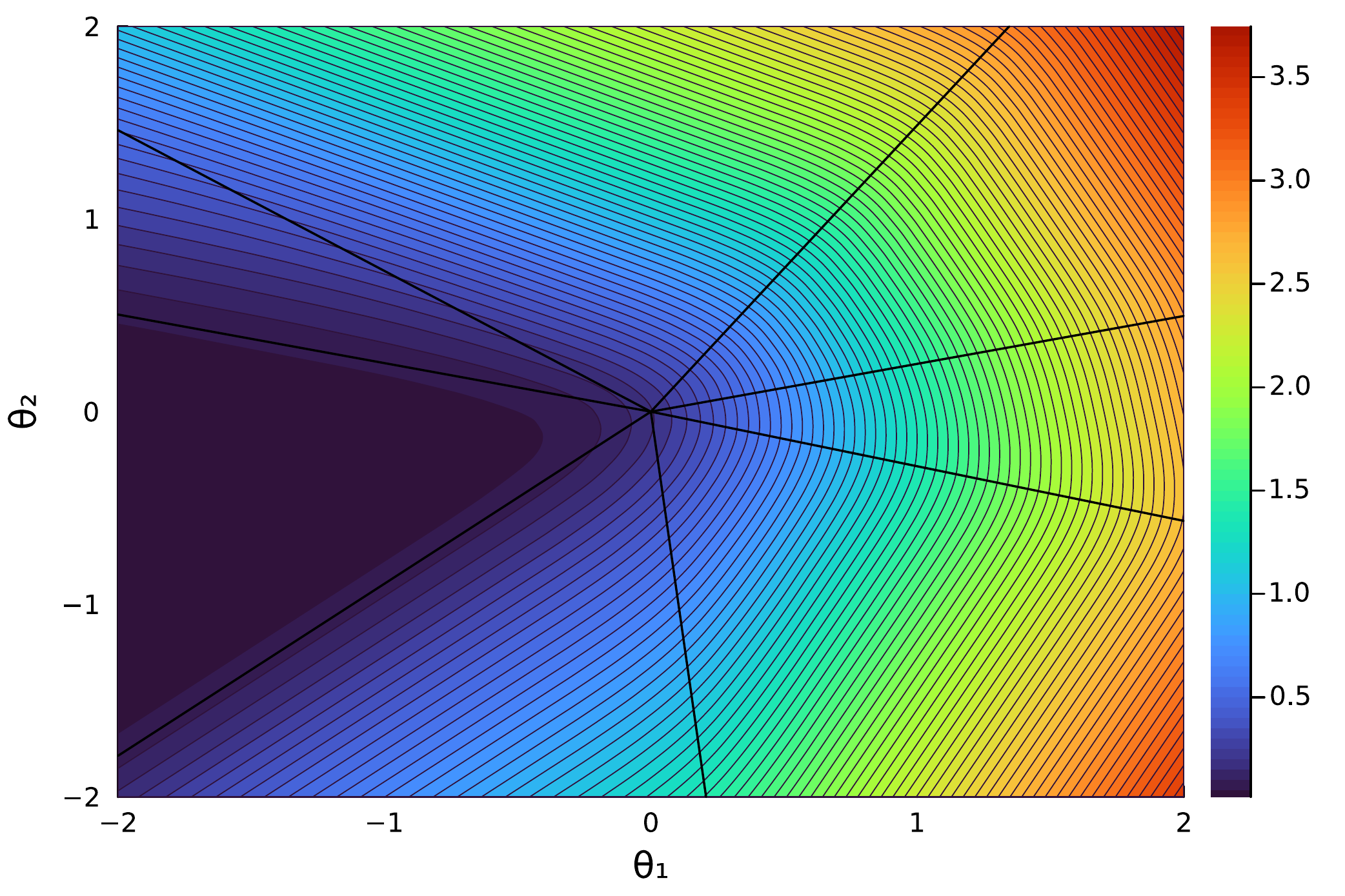}
        \caption{$\theta\mapsto\calL_\varepsilon(\theta, \bar y)$ for $h=0$}
        \label{subfig:FYL_perturbed}
    \end{subfigure}
    \begin{subfigure}[b]{0.45\textwidth}    
        \includegraphics[width=\textwidth]{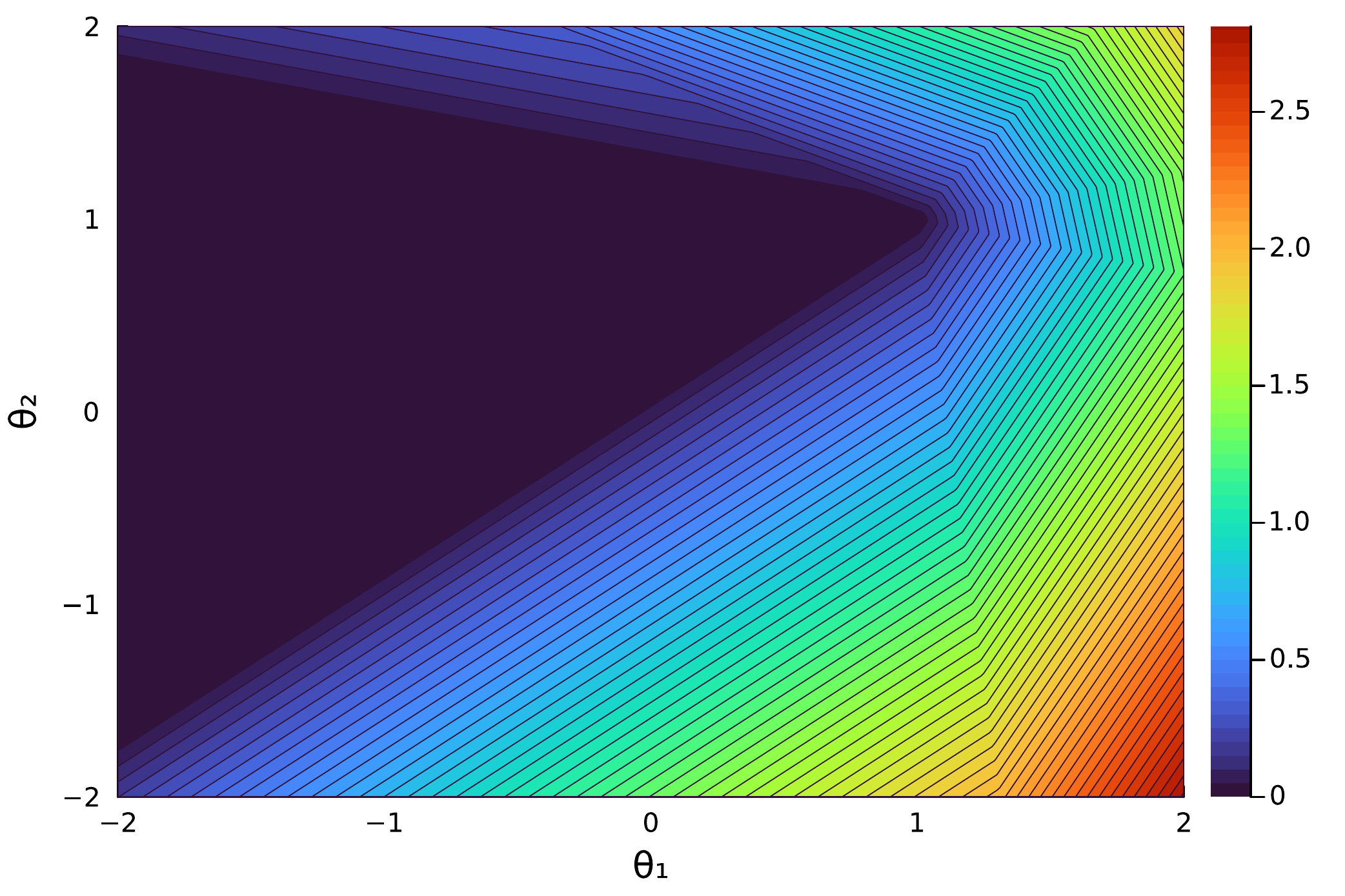}
        \caption{$\theta\mapsto\calL(\theta,\bar y)$ for $h\neq0$}
        \label{subfig:FYL_d}
    \end{subfigure}
    \begin{subfigure}[b]{0.45\textwidth}    
        \includegraphics[width=\textwidth]{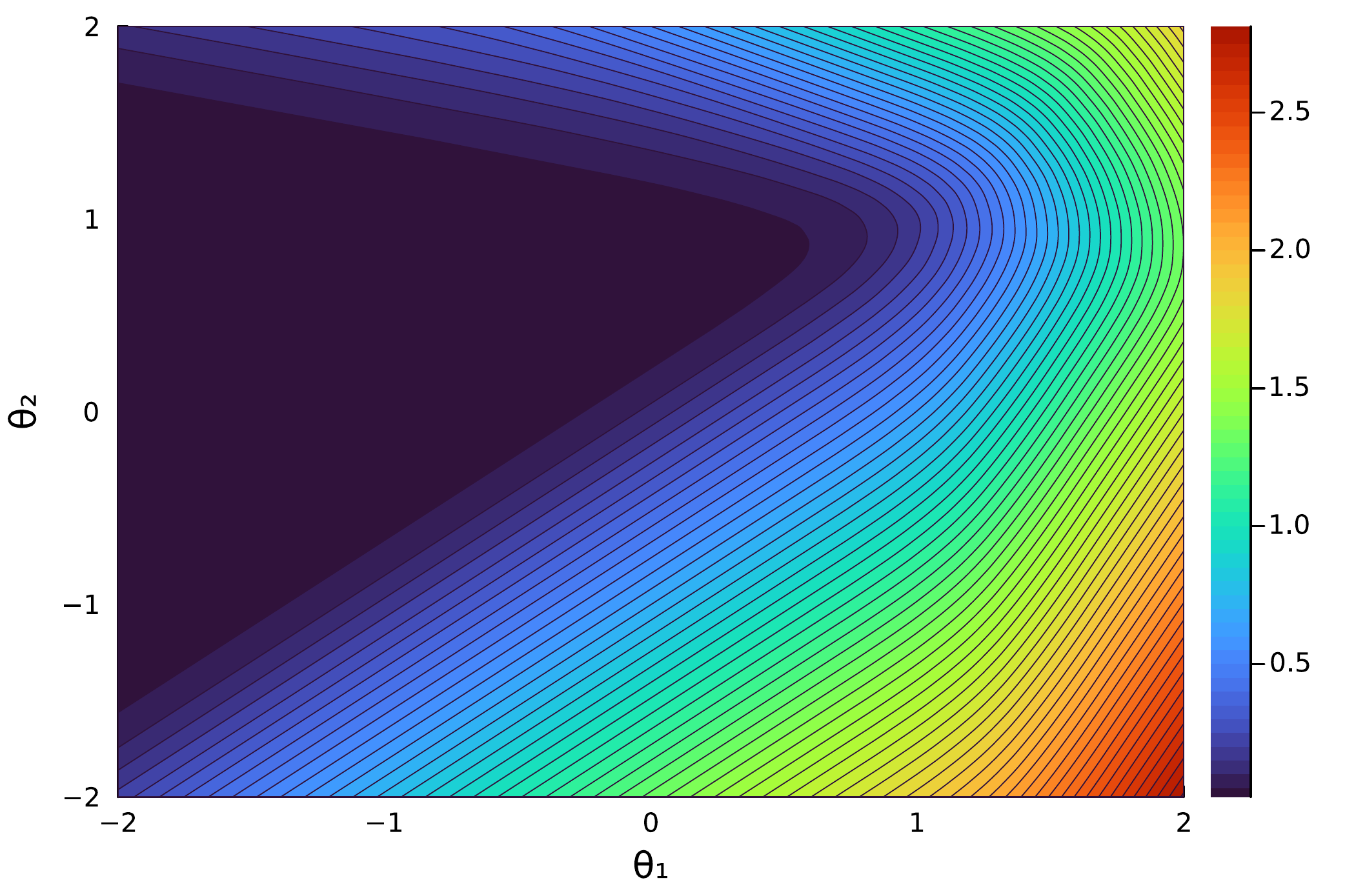}
        \caption{$\theta\mapsto\calL_\varepsilon(\theta,\bar y)$ for $h\neq0$}
        \label{subfig:FYL_e}
    \end{subfigure}
    \caption{Example with two-dimensional $\theta$ and its polytope~\ref{subfig:polytope}. Contour plots of the loss $\calL$ and its perturbed version $\calL$ for two different values of $h$. Thick lines on~\ref{subfig:FYL} and~\ref{subfig:FYL_perturbed} represent the normal fan associated to the polytope. We can see the perturbed regularization ``pushing'' the loss inside the normal cone $\calC(\bar y)$.}
    \label{fig:2D_FYL}
\end{figure}

The perturbed loss $\mathcal{L}_\varepsilon$ has properties that make it very suitable for a learning problem.
Indeed, it is convex and smooth, and while the expectation in its gradient is intractable, sampling $Z$ gives a stochastic gradient. 
The learning problem~\eqref{eq:learning-problem} with loss $\mathcal{L}_\varepsilon$ can therefore be solved using stochastic gradient descent.
Furthermore, $\theta\mapsto\calL_\varepsilon(\theta, \bar y)$ tends to $0$ only if $\theta$ is far in the interior of $\mathcal{C}(\bar y)$, which means that there is no ambiguity on the fact that $\bar y$ is an optimal solution of~\eqref{eq:co_layer} for $\theta$, and removes the degeneracy issue of $\theta\mapsto\calL(\theta,\bar y)$.
This situation is illustrated in two dimensions in Figure~\ref{fig:2D_FYL}.
Figure~\ref{subfig:polytope} represents the polytope corresponding to the convex envelope $\conv(g(\calY))) = \conv(\{g(y),\, \forall y\in\calY\})$. 
The algorithm~$f$ necessarily outputs a vertex of this polytope (red square) as the objective function of the \gls*{CO}-layer is linear in~$g(y)$.
The dark blue hexagon is the output $\mathbb{E}[g(f(\theta+\varepsilon Z))]$ of the perturbed \gls*{CO}-layer, which can be seen as the expectation over a distribution (blue circles) on the vertices of the polytope.
On Figure~\ref{subfig:FYL}, with $h=0$ and no perturbation, we can see that the loss $\calL(\theta,\bar y)$ is $0$ for any $\theta$ in the normal cone of $\bar y$ (with $g(\bar y)$ being the red square), and that the origin belongs to the normal cone of any vertex $y$.
On Figure~\ref{subfig:FYL_perturbed}, the perturbation has been added and we can see that $\calL(\theta,\bar y)$ goes to zero only when we are safely inside the cone, which means that the output of the combinatorial optimization layer is non-ambiguous. Finally, on Figure~\ref{subfig:FYL_d} and \ref{subfig:FYL_e}, we can see that a non-zero $h$ changes the geometry for small $\theta$ values but not for large $\theta$ values.

\paragraph{Non-optimal CO-layer} Note that our \gls*{CO}-layer $\COptimizer$ is a (meta)-heuristic algorithm and thus does not guarantee an optimal solution, which is not in line with the assumptions of the theory reviewed in this section. However, we observe that this is not a problem in practice as $\COptimizer$ outputs solutions close enough to optimal ones.

\section{Computational study}\label{sec:results}

The aim of our computational study is twofold. First, we validate the performance of our \gls*{CO}-enriched \gls*{ML} pipeline in a benchmark against several state-of-the-art approaches. Second, we conduct extensive numerical experiments to assess the impact of different training settings on the performance of our \gls*{CO}-enriched \gls*{ML} pipeline. Specifically, we investigate the impact of i) the feature set, ii) the size of instances in the training set, iii) the training set size, iv) the imitated target strategy, and v) the type of statistical model used.

For this purpose, we first detail the design of our computational study in Section~\ref{subsec:design_of_experiments}, and the different benchmarks in Section~\ref{subsec:benchmarks}. We then present the results of our benchmark study in Section~\ref{subsec:analysis}, and discuss the results of our experiments on different training settings in Section~\ref{subsec:extendedAnalysis}.

\subsection{Design of experiments}\label{subsec:design_of_experiments}
We design our computational study in line with the problem setting presented in the EURO Meets NeurIPS Vehicle Routing Competition \parencite{euromeetsneurips2022}.
We consider a set of \gls*{VRPTW} instances derived from real-world data of a US-based grocery delivery service. We refer to these instances as \emph{static} instances. A static instance contains a set of requests, each with a specific location, a service time, a demand, and a time window. Hence, each instance implicitly defines a distribution of request locations, service times, demands, and time windows. We generate an instance of the \emph{dynamic} problem from a static instance as follows:
\noindent we first discretize the planning horizon of the static instance into one hour epochs. Then, for each epoch $\Epoch$, we sample $m$ requests randomly, constructed by drawing from the static instance's location, demand, service time, and time window distributions independently. We refer to $m$ as the \emph{sample size} in the remainder of this section.
We discard requests that are infeasible with respect to the starting time of the current epoch and the chosen time window. The remaining requests form the set of requests revealed in epoch~$\Epoch$. This sampling approach has two implications: First, the expected number of requests arriving in each epoch decreases with advancing epochs as the probability of sampling a feasible time window decreases. Second, time windows that end closer to the end of the planning horizon are more likely to appear, as these time windows have a higher probability to be feasible.
In line with the challenge's problem setting, we minimize travel time only, such that the number of vehicles, service times, and waiting times are not part of the objective. Furthermore, we do not limit the number of vehicles and assume knowledge over the static instance's request distribution.

Finally, we note that we can generate several dynamic instances from one static instance by varying the seed used to sample from the request distribution. In what follows, we refer to this seed as the \emph{instance seed}.

\subsection{Benchmark policies}\label{subsec:benchmarks}

We evaluate the performance of our \gls*{CO}-enriched \gls*{ML} pipeline against a set of benchmark policies. These policies follow a two-stage approach by first deciding on the set of requests to dispatch and then form covering routes in a second step by solving a \gls*{VRPTW} on the dispatched requests using the \gls*{HGS} algorithm \parencite[cf.][]{koolHybrid2022, vidalHybrid2021}.
For each benchmark policy, we allow a time limit for finding decisions of 90 seconds per epoch, if not stated differently.
In what follows, we briefly introduce the benchmark policies considered in this computational study.

\paragraph{Greedy policy}
The greedy policy dispatches all requests as soon as they enter the system.

\paragraph{Lazy policy}
The greedy policy dispatches all requests as soon as they enter the system.

\paragraph{Random policy}
The random policy dispatches postponable requests with a probability of~50\%. It further dispatches all must-dispatch requests.

\paragraph{Rolling-horizon policy}
The rolling-horizon policy samples a scenario for the remaining epochs, and applies the \gls*{HGS} in a similar way as the anticipative strategy assuming that the sampled scenario represents the true scenario. This yields a set of routes based on which it decides which requests to dispatch. Specifically, it dispatches a request of the current epoch if and only if that request shares a route with a must dispatch request of the current epoch. We assign a time limit of 600 seconds to ensure convergence of the \gls*{HGS} since this policy entails solving a \gls*{VRPTW} on a complete scenario. Note that this extends the 90 second time

\paragraph{Monte-carlo policy}
The monte-carlo policy samples nine scenarios for the remaining epochs, solved individually as in the rolling-horizon policy. It dispatches a request based on a majority decision, i.e., if and only if the rolling-horizon policy dispatches the request in at least 50\% of the sampled scenarios. We raise the time-limit accordingly, i.e., allow a total of 5400 seconds.

\paragraph{ML-CO policy}
The ML-CO policy is encoded in the CO-enriched ML pipeline introduced in Section \ref{sec:pipeline}. 
Unless specified otherwise, we train our ML-CO policy on a set of 15 training instances using a sample size of 50 requests per epoch, derived using the anticipative strategy as detailed in Section~\ref{subsec:anticipative_policy}. We limit the runtime of the \gls*{HGS} used to derive the anticipative solutions to 3600 seconds.
Our \gls*{ML}-layer uses a feedforward neural network with four hidden layers, each comprising ten neurons.
For the ML-CO policy, we set the time limit for finding decisions to 90 seconds per epoch. Table~\ref{tab:runtime} indicates that this time limit does not impact the performance of our ML-CO policy significantly.
\begin{table}[h]
    \footnotesize
    \centering
    \begin{tabular}{c | c c c c c c} 
     & \multicolumn{6}{c}{Decision time limit during evaluation [seconds]} \\
     & 30 & 60 & 90 & 120 & 180 & 240 \\ 
     \hline\hline
     \textbf{Relative distance to ant. b.:} & 5.66\% & 5.18\% & 5.15\% & 5.01\% & 4.89\% & 4.87\% \\
     \hline
    \end{tabular}
    \caption{Performance of our ML-CO policy for different time limits that we allow for the PC-HGS during evaluation.}
    \label{tab:runtime}
\end{table}

\bigskip
We note that the \textit{greedy}, \textit{lazy}, and \textit{random} policies are the baseline heuristics used in the EURO Meets NeurIPS Vehicle Routing Competition and refer the interested reader to \cite{euromeetsneurips2022} for more details.

We additionally include the anticipative strategy used during training (cf. Section~\ref{subsec:anticipative_policy}) as a baseline. To derive the anticipative baseline, we solve the offline \gls*{VRPTW} using the \gls*{HGS} algorithm proposed in \cite{koolHybrid2022} with a time limit of 3600 seconds. Note that the \gls*{HGS} algorithm is a heuristic which implies that 
it is possible to find a solution which outperforms the anticipative baseline.

We evaluate the performance of the benchmark policies on~25 dynamic test instances and report for each instance the average result over 20 different instance seeds. To create a dynamic test instance we consider a sample size of 100 requests per epoch.

\subsection{Performance Analysis}
\label{subsec:analysis}

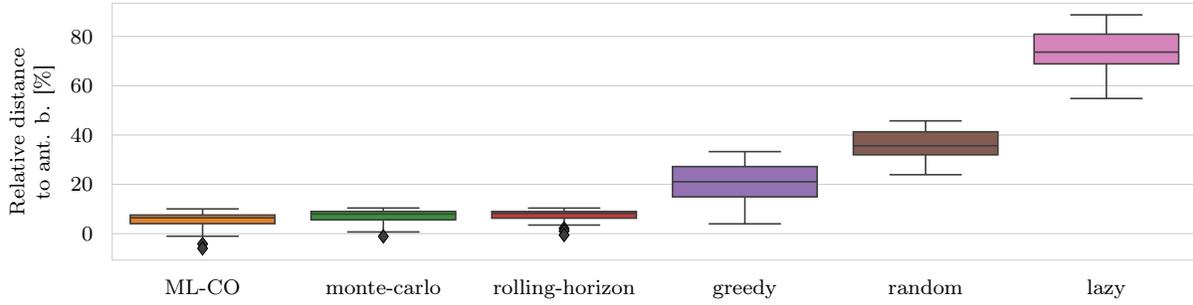
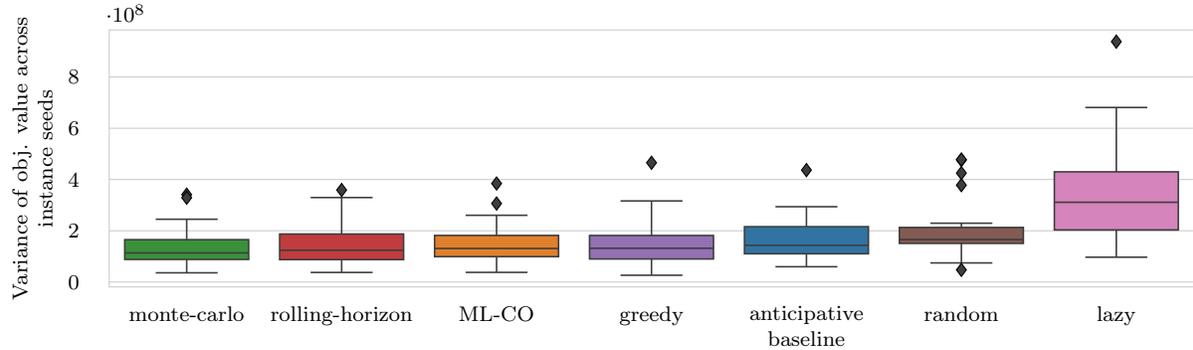
\begin{figure}[ht]
    \centering
    \footnotesize
    \begin{subfigure}[b]{\textwidth}
         \centering
\begin{tikzpicture}

\definecolor{brown1926061}{RGB}{192,60,61}
\definecolor{darkslategray38}{RGB}{38,38,38}
\definecolor{darkslategray61}{RGB}{61,61,61}
\definecolor{dimgray1319183}{RGB}{131,91,83}
\definecolor{lightgray204}{RGB}{204,204,204}
\definecolor{mediumpurple147113178}{RGB}{147,113,178}
\definecolor{orchid213132188}{RGB}{213,132,188}
\definecolor{peru22412844}{RGB}{224,128,44}
\definecolor{seagreen5814558}{RGB}{58,145,58}

\begin{axis}[
axis line style={lightgray204},
height=5cm,
scaled y ticks=true,
tick align=outside,
width=16cm,
x grid style={lightgray204},
xmajorticks=false,
xmajorticks=true,
xmin=-0.5, xmax=5.5,
xtick style={color=darkslategray38},
xtick style={draw=none},
xtick={0,1,2,3,4,5},
xticklabel style={align=center},
xticklabels={ML-CO,monte-carlo,rolling-horizon,greedy,random,lazy},
y grid style={lightgray204},
ylabel style={align=center},
ylabel={Relative distance \\ to ant. b. [\%]},
ymajorgrids,
ymajorticks=false,
ymajorticks=true,
ymin=-10.7068181156302, ymax=93.5083887724209,
ytick style={color=darkslategray38},
ytick style={draw=none}
]
\path [draw=darkslategray61, fill=peru22412844, semithick]
(axis cs:-0.4,4.07950156378079)
--(axis cs:0.4,4.07950156378079)
--(axis cs:0.4,7.57757278779952)
--(axis cs:-0.4,7.57757278779952)
--(axis cs:-0.4,4.07950156378079)
--cycle;
\path [draw=darkslategray61, fill=seagreen5814558, semithick]
(axis cs:0.6,5.62378430399234)
--(axis cs:1.4,5.62378430399234)
--(axis cs:1.4,9.03532422863164)
--(axis cs:0.6,9.03532422863164)
--(axis cs:0.6,5.62378430399234)
--cycle;
\path [draw=darkslategray61, fill=brown1926061, semithick]
(axis cs:1.6,6.29875810608642)
--(axis cs:2.4,6.29875810608642)
--(axis cs:2.4,9.00905510027842)
--(axis cs:1.6,9.00905510027842)
--(axis cs:1.6,6.29875810608642)
--cycle;
\path [draw=darkslategray61, fill=mediumpurple147113178, semithick]
(axis cs:2.6,14.9434699868178)
--(axis cs:3.4,14.9434699868178)
--(axis cs:3.4,27.2316924282875)
--(axis cs:2.6,27.2316924282875)
--(axis cs:2.6,14.9434699868178)
--cycle;
\path [draw=darkslategray61, fill=dimgray1319183, semithick]
(axis cs:3.6,31.995954025116)
--(axis cs:4.4,31.995954025116)
--(axis cs:4.4,41.3380706478463)
--(axis cs:3.6,41.3380706478463)
--(axis cs:3.6,31.995954025116)
--cycle;
\path [draw=darkslategray61, fill=orchid213132188, semithick]
(axis cs:4.6,68.8895452918279)
--(axis cs:5.4,68.8895452918279)
--(axis cs:5.4,80.945837838266)
--(axis cs:4.6,80.945837838266)
--(axis cs:4.6,68.8895452918279)
--cycle;
\addplot [semithick, darkslategray61]
table {%
0 4.07950156378079
0 -1.01029782551158
};
\addplot [semithick, darkslategray61]
table {%
0 7.57757278779952
0 10.0500268225354
};
\addplot [semithick, darkslategray61]
table {%
-0.2 -1.01029782551158
0.2 -1.01029782551158
};
\addplot [semithick, darkslategray61]
table {%
-0.2 10.0500268225354
0.2 10.0500268225354
};
\addplot [black, mark=diamond*, mark size=2.5, mark options={solid,fill=darkslategray61}, only marks]
table {%
0 -4.14832007418244
0 -5.96976325708246
};
\addplot [semithick, darkslategray61]
table {%
1 5.62378430399234
1 0.702961017016277
};
\addplot [semithick, darkslategray61]
table {%
1 9.03532422863164
1 10.4141846508041
};
\addplot [semithick, darkslategray61]
table {%
0.8 0.702961017016277
1.2 0.702961017016277
};
\addplot [semithick, darkslategray61]
table {%
0.8 10.4141846508041
1.2 10.4141846508041
};
\addplot [black, mark=diamond*, mark size=2.5, mark options={solid,fill=darkslategray61}, only marks]
table {%
1 -1.09858197075344
};
\addplot [semithick, darkslategray61]
table {%
2 6.29875810608642
2 3.50306332208783
};
\addplot [semithick, darkslategray61]
table {%
2 9.00905510027842
2 10.3928104429071
};
\addplot [semithick, darkslategray61]
table {%
1.8 3.50306332208783
2.2 3.50306332208783
};
\addplot [semithick, darkslategray61]
table {%
1.8 10.3928104429071
2.2 10.3928104429071
};
\addplot [black, mark=diamond*, mark size=2.5, mark options={solid,fill=darkslategray61}, only marks]
table {%
2 2.14608842826696
2 1.1333983789793
2 -0.543651026411002
};
\addplot [semithick, darkslategray61]
table {%
3 14.9434699868178
3 3.99819767420397
};
\addplot [semithick, darkslategray61]
table {%
3 27.2316924282875
3 33.2888801053517
};
\addplot [semithick, darkslategray61]
table {%
2.8 3.99819767420397
3.2 3.99819767420397
};
\addplot [semithick, darkslategray61]
table {%
2.8 33.2888801053517
3.2 33.2888801053517
};
\addplot [semithick, darkslategray61]
table {%
4 31.995954025116
4 23.9275976933183
};
\addplot [semithick, darkslategray61]
table {%
4 41.3380706478463
4 45.7699430448815
};
\addplot [semithick, darkslategray61]
table {%
3.8 23.9275976933183
4.2 23.9275976933183
};
\addplot [semithick, darkslategray61]
table {%
3.8 45.7699430448815
4.2 45.7699430448815
};
\addplot [semithick, darkslategray61]
table {%
5 68.8895452918279
5 54.8517223415442
};
\addplot [semithick, darkslategray61]
table {%
5 80.945837838266
5 88.7713339138731
};
\addplot [semithick, darkslategray61]
table {%
4.8 54.8517223415442
5.2 54.8517223415442
};
\addplot [semithick, darkslategray61]
table {%
4.8 88.7713339138731
5.2 88.7713339138731
};
\addplot [semithick, darkslategray61]
table {%
-0.4 6.4643243249672
0.4 6.4643243249672
};
\addplot [semithick, darkslategray61]
table {%
0.6 7.95750119666611
1.4 7.95750119666611
};
\addplot [semithick, darkslategray61]
table {%
1.6 8.23894086719628
2.4 8.23894086719628
};
\addplot [semithick, darkslategray61]
table {%
2.6 21.0639969673116
3.4 21.0639969673116
};
\addplot [semithick, darkslategray61]
table {%
3.6 35.681984875208
4.4 35.681984875208
};
\addplot [semithick, darkslategray61]
table {%
4.6 73.6540750667643
5.4 73.6540750667643
};
\end{axis}

\end{tikzpicture}
         \caption{Performance relative to anticipative baseline solution.}
         \label{fig:resultsRelative}
     \end{subfigure}
     \hfill
     \begin{subfigure}[b]{\textwidth}
         \centering
\begin{tikzpicture}

\definecolor{brown1926061}{RGB}{192,60,61}
\definecolor{darkslategray38}{RGB}{38,38,38}
\definecolor{darkslategray61}{RGB}{61,61,61}
\definecolor{dimgray1319183}{RGB}{131,91,83}
\definecolor{lightgray204}{RGB}{204,204,204}
\definecolor{mediumpurple147113178}{RGB}{147,113,178}
\definecolor{orchid213132188}{RGB}{213,132,188}
\definecolor{peru22412844}{RGB}{224,128,44}
\definecolor{seagreen5814558}{RGB}{58,145,58}
\definecolor{steelblue49115161}{RGB}{49,115,161}

\begin{axis}[
axis line style={lightgray204},
height=5cm,
scaled y ticks=true,
tick align=outside,
width=16cm,
x grid style={lightgray204},
xmajorticks=false,
xmajorticks=true,
xmin=-0.5, xmax=6.5,
xtick style={color=darkslategray38},
xtick style={draw=none},
xtick={0,1,2,3,4,5,6},
xticklabel style={align=center},
xticklabels={monte-carlo,rolling-horizon,ML-CO,greedy,anticipative \\ baseline,random,lazy},
y grid style={lightgray204},
ylabel style={align=center},
ylabel={Variance of obj. value across \\ instance seeds},
ymajorgrids,
ymajorticks=false,
ymajorticks=true,
ymin=-19062785.986, ymax=983347714.886,
ytick style={color=darkslategray38},
ytick style={draw=none},
ytick={-200000000,0,200000000,400000000,600000000,800000000,1000000000},
yticklabel style={xshift=-4pt},
yticklabels={\ensuremath{-}2,0,2,4,6,8,10}
]
\path [draw=darkslategray61, fill=seagreen5814558, semithick]
(axis cs:-0.4,88118614.19)
--(axis cs:0.4,88118614.19)
--(axis cs:0.4,165268941.8275)
--(axis cs:-0.4,165268941.8275)
--(axis cs:-0.4,88118614.19)
--cycle;
\path [draw=darkslategray61, fill=brown1926061, semithick]
(axis cs:0.6,87379179.7475)
--(axis cs:1.4,87379179.7475)
--(axis cs:1.4,187391569.3475)
--(axis cs:0.6,187391569.3475)
--(axis cs:0.6,87379179.7475)
--cycle;
\path [draw=darkslategray61, fill=peru22412844, semithick]
(axis cs:1.6,98952432.9275)
--(axis cs:2.4,98952432.9275)
--(axis cs:2.4,181779578.16)
--(axis cs:1.6,181779578.16)
--(axis cs:1.6,98952432.9275)
--cycle;
\path [draw=darkslategray61, fill=mediumpurple147113178, semithick]
(axis cs:2.6,89900861.79)
--(axis cs:3.4,89900861.79)
--(axis cs:3.4,181588032.0275)
--(axis cs:2.6,181588032.0275)
--(axis cs:2.6,89900861.79)
--cycle;
\path [draw=darkslategray61, fill=steelblue49115161, semithick]
(axis cs:3.6,110212568.71)
--(axis cs:4.4,110212568.71)
--(axis cs:4.4,215937517.8475)
--(axis cs:3.6,215937517.8475)
--(axis cs:3.6,110212568.71)
--cycle;
\path [draw=darkslategray61, fill=dimgray1319183, semithick]
(axis cs:4.6,150673466.36)
--(axis cs:5.4,150673466.36)
--(axis cs:5.4,212926818.3275)
--(axis cs:4.6,212926818.3275)
--(axis cs:4.6,150673466.36)
--cycle;
\path [draw=darkslategray61, fill=orchid213132188, semithick]
(axis cs:5.6,203220172.66)
--(axis cs:6.4,203220172.66)
--(axis cs:6.4,429888234.04)
--(axis cs:5.6,429888234.04)
--(axis cs:5.6,203220172.66)
--cycle;
\addplot [semithick, darkslategray61]
table {%
0 88118614.19
0 35978580.4275
};
\addplot [semithick, darkslategray61]
table {%
0 165268941.8275
0 244991558.4875
};
\addplot [semithick, darkslategray61]
table {%
-0.2 35978580.4275
0.2 35978580.4275
};
\addplot [semithick, darkslategray61]
table {%
-0.2 244991558.4875
0.2 244991558.4875
};
\addplot [black, mark=diamond*, mark size=2.5, mark options={solid,fill=darkslategray61}, only marks]
table {%
0 340902155.29
0 328885739.4475
};
\addplot [semithick, darkslategray61]
table {%
1 87379179.7475
1 37388691.96
};
\addplot [semithick, darkslategray61]
table {%
1 187391569.3475
1 329563114.0275
};
\addplot [semithick, darkslategray61]
table {%
0.8 37388691.96
1.2 37388691.96
};
\addplot [semithick, darkslategray61]
table {%
0.8 329563114.0275
1.2 329563114.0275
};
\addplot [black, mark=diamond*, mark size=2.5, mark options={solid,fill=darkslategray61}, only marks]
table {%
1 358983427.5475
};
\addplot [semithick, darkslategray61]
table {%
2 98952432.9275
2 37667208.05
};
\addplot [semithick, darkslategray61]
table {%
2 181779578.16
2 260401197.1275
};
\addplot [semithick, darkslategray61]
table {%
1.8 37667208.05
2.2 37667208.05
};
\addplot [semithick, darkslategray61]
table {%
1.8 260401197.1275
2.2 260401197.1275
};
\addplot [black, mark=diamond*, mark size=2.5, mark options={solid,fill=darkslategray61}, only marks]
table {%
2 306290135.0275
2 384334814.6475
};
\addplot [semithick, darkslategray61]
table {%
3 89900861.79
3 26501327.69
};
\addplot [semithick, darkslategray61]
table {%
3 181588032.0275
3 316269068.96
};
\addplot [semithick, darkslategray61]
table {%
2.8 26501327.69
3.2 26501327.69
};
\addplot [semithick, darkslategray61]
table {%
2.8 316269068.96
3.2 316269068.96
};
\addplot [black, mark=diamond*, mark size=2.5, mark options={solid,fill=darkslategray61}, only marks]
table {%
3 465330189.54
};
\addplot [semithick, darkslategray61]
table {%
4 110212568.71
4 59673177.2475
};
\addplot [semithick, darkslategray61]
table {%
4 215937517.8475
4 293827195.56
};
\addplot [semithick, darkslategray61]
table {%
3.8 59673177.2475
4.2 59673177.2475
};
\addplot [semithick, darkslategray61]
table {%
3.8 293827195.56
4.2 293827195.56
};
\addplot [black, mark=diamond*, mark size=2.5, mark options={solid,fill=darkslategray61}, only marks]
table {%
4 436658588.3275
};
\addplot [semithick, darkslategray61]
table {%
5 150673466.36
5 74566578.5275
};
\addplot [semithick, darkslategray61]
table {%
5 212926818.3275
5 229573082.7475
};
\addplot [semithick, darkslategray61]
table {%
4.8 74566578.5275
5.2 74566578.5275
};
\addplot [semithick, darkslategray61]
table {%
4.8 229573082.7475
5.2 229573082.7475
};
\addplot [black, mark=diamond*, mark size=2.5, mark options={solid,fill=darkslategray61}, only marks]
table {%
5 47359593.7475
5 475981006.8475
5 377868749.31
5 425026322.56
5 478049704.94
};
\addplot [semithick, darkslategray61]
table {%
6 203220172.66
6 96933985.19
};
\addplot [semithick, darkslategray61]
table {%
6 429888234.04
6 680973012.6875
};
\addplot [semithick, darkslategray61]
table {%
5.8 96933985.19
6.2 96933985.19
};
\addplot [semithick, darkslategray61]
table {%
5.8 680973012.6875
6.2 680973012.6875
};
\addplot [black, mark=diamond*, mark size=2.5, mark options={solid,fill=darkslategray61}, only marks]
table {%
6 937783601.21
};
\addplot [semithick, darkslategray61]
table {%
-0.4 113735868.0875
0.4 113735868.0875
};
\addplot [semithick, darkslategray61]
table {%
0.6 123886880.4275
1.4 123886880.4275
};
\addplot [semithick, darkslategray61]
table {%
1.6 130985973.8275
2.4 130985973.8275
};
\addplot [semithick, darkslategray61]
table {%
2.6 131764172.5
3.4 131764172.5
};
\addplot [semithick, darkslategray61]
table {%
3.6 143281447.7275
4.4 143281447.7275
};
\addplot [semithick, darkslategray61]
table {%
4.6 165684652.84
5.4 165684652.84
};
\addplot [semithick, darkslategray61]
table {%
5.6 311174306.1275
6.4 311174306.1275
};
\end{axis}

\end{tikzpicture}
         \caption{Variance of objective value with respect to instance seeds.}
         \label{fig:resultsVar}
     \end{subfigure}
    \caption{Performance of benchmark policies.}
\end{figure}

Figure \ref{fig:resultsRelative} compares the performances of the introduced benchmarks. The figure shows the gap of the policy's objective value relative to the anticipative baseline~(ant. b.). 
In general, we can see that benchmark policies which consider information about the uncertain appearance of future requests, i.e., rolling-horizon, monte-carlo, ML-CO, and the anticipative baseline, outperform benchmarks which make dispatching decision based on the current epoch only, i.e., lazy, random, and greedy.
Specifically, these perform on average~74.05\%,~36.29\%, and~20.99\% worse than the anticipative baseline.
The \textit{rolling-horizon} policy performs 7.12\% worse than the anticipative baseline. We see that considering the future impact of our dispatching decision based on only a single scenario already improves the performance in comparison to the greedy policy significantly.
This is surprising as the scenario drawn is often only a poor representation of the requests actually observed in later epochs due to the size of the sample space. However, the distance between the drawn scenario and the correct scenario might be outweighed by the benefit of combining future requests with actual requests in low-cost routes.
The \textit{monte-carlo} policy outperforms the rolling-horizon policy and only has a gap of 6.97\% to the anticipative baseline. This shows that a better approximation of future uncertainties, achieved through a higher number of sampled scenarios, leads to better dispatching decisions, thus improving performance. 
Yet, the improvement with respect to the rolling-horizon benchmark is rather small in comparison to the improvement from the greedy policy to the rolling-horizon policy.
Our \textit{ML-CO} policy outperforms all other online policies and is only 5.15\% worse than the anticipative baseline.
This is surprising as the general consensus in literature on multi-stage stochastic problems indicates the contrary, i.e., that imitating an anticipative strategy does not generalize well.
\begin{result}
Our ML-CO policy performs best across all online benchmark policies and outperforms the monte-carlo policy by 1.57\%. This indicates that learning a policy by imitating an anticipative strategy yields good performances in this problem setting.
\end{result}

Figure \ref{fig:resultsVar} shows the variance of the objective value over several instance seeds for the considered benchmark policies. Comparing Figures~\ref{fig:resultsRelative} and~\ref{fig:resultsVar} shows a clear trend: policies with low variance outperform policies with high variance.
This trend results from finding more robust solutions that generalize well over uncertain future observations. 
Note, that, although performing worse, the monte-carlo policy has less variance than our ML-CO policy.
This relates to the sampling approach the monte-carlo policy bases on, which yields a robust solution that performs well over all sampled scenarios and therefore focuses too much on variance minimization. Our ML-CO policy on the other hand balances the trade-off between robustness and solution quality.

\begin{result}
The monte-carlo policy and our ML-CO policy generalize well over uncertain future observations. The learning component in the ML-CO policy balances the trade-off between robustness and solution quality, leading to a superior performance.
\end{result}

\subsection{Extended analysis}
\label{subsec:extendedAnalysis}
This analysis assesses the impact of different training settings on the performance of our ML-CO policy. Specifically, we investigate the impact of i) the feature set, ii) the size of each training instance, iii) the training set size, iv) the imitated target strategy, and v) the type of statistical model used.
We evaluate each model's performance as detailed in Section~\ref{subsec:benchmarks} and report the relative gap to the anticipative baseline.

\paragraph{Different feature sets}
\begin{table}[ht]
\footnotesize
\centering
\begin{tabular}{c | c c c} 
 & \multicolumn{3}{c}{Feature sets} \\
 & complete & model-aware & model-free \\ 
 & (baseline) & & \\ 
 \hline\hline
 \textbf{Relative distance to ant. b.:} & 5.15\% & 13.83\% & 6.78\% \\
 \hline
\end{tabular}
\caption{Performance of ML-CO policy using different feature sets.}
\label{tab:featureSetPerformance}
\end{table}
Table \ref{tab:featureSetPerformance} compares the performance of our ML-CO policy on three feature sets (i.e., \emph{complete}, \emph{model-aware}, and \emph{model-free}).
Table~\ref{tab:featureSets} details the features each set comprises.
Specifically, the \emph{model-free} feature set contains features computable from the current state $x^e$, while the \emph{model-aware} feature set only contains features which include distributional information from the static instance. We further include the \emph{complete} feature set which combines the features from the \emph{model-aware} and \emph{model-free} feature sets. The \emph{complete} feature set was used in the model submitted to the challenge, i.e., our baseline.
Our results show that the performance of our ML-CO policy does not rely solely on \emph{model-aware} information not available in real-world scenarios. Specifically, considering \emph{model-aware} features derived from the static instance decreases the gap of our ML-CO policy to the anticipative baseline by only $1.63$ percentage points on average.
\begin{result}
Our ML-CO policy performs well without considering \emph{model-aware} features derived from the static instance.
\end{result}
\begin{table}[ht]
\footnotesize
\centering
\begin{tabular}{l l | l l} 
 \multicolumn{2}{c}{model-free} & \multicolumn{2}{c}{model-aware} \\
 \hline\hline
    x coordinate of location & $x_r$ & \multicolumn{2}{l}{\textit{Quantiles from distribution of travel time to all locations:}} \\
    y coordinate of location & $y_r$ & 1\% quantile & $Pr[X<x] \leq 0.01, X\sim t_{r,:}$ \\
    demand & $q_r$ & 5\% quantile & $Pr[X<x] \leq 0.05, X\sim t_{r,:}$ \\
    service time & $s_r$ & 10\% quantile & $Pr[X<x] \leq 0.1, X\sim t_{r,:}$ \\
    time window start & $l_r$ & 50\% quantile & $Pr[X<x] \leq 0.5, X\sim t_{r,:}$ \\
    time window end & $u_r$ & \multicolumn{2}{l}{\textit{Quantiles from distribution of slack time to all time windows:}} \\
    time from depot to request & $t_{d, r}$ & 0\% quantile & $Pr[X<x] \leq 0, X\sim u_: - (l_r+s_r + t_{r,:})$ \\
    relative time depot to request &$t_{d, r} / (u_r - s_r)$ & 1\% quantile & $Pr[X<x] \leq 0.01, X\sim u_: - (l_r+s_r + t_{r,:})$ \\
    time window start / rem. time & $l_r / (T_{max} - \TimeStep_{\Epoch})$ & 5\% quantile & $Pr[X<x] \leq 0.05, X\sim u_: - (l_r+s_r+ t_{r,:})$ \\
    time window end / rem. time & $u_r / (T_{max} - \TimeStep_{\Epoch})$ & 10\% quantile & $Pr[X<x] \leq 0.1, X\sim u_: - (l_r+s_r+ t_{r,:})$ \\
    is must dispatch & $\mathbbm{1}_{\TimeStep_{\Epoch}+\DispatchTime + t_{d, r}>u_r}$ & 50\% quantile & $Pr[X<x] \leq 0.5, X\sim u_: - (l_r+s_r+ t_{r,:})$ \\
 \hline
\end{tabular}
\caption{Different feature sets.}
\label{tab:featureSets}
\end{table}

\paragraph{Different sample size of training instances}
\begin{table}[t]
\footnotesize
\centering
\begin{tabular}{c | c c c c c} 
 & \multicolumn{5}{c}{Sample size of training instances} \\
 & 10 & 25 & 50 & 75 & 100 \\ 
 \hline\hline
 \textbf{Relative distance to ant. b.:} & 8.71\% & 6.31\% & 5.15\% & 7.41\% & 13.08\% \\
 \hline
\end{tabular}
\caption{Performance of ML-CO policy when trained on different sized training instances.}
\label{tab:instanceSizes}
\end{table}
Table~\ref{tab:instanceSizes} shows the performance of the ML-CO policy relative to the performance of the anticipative baseline when training the ML-CO policy on training instances derived from different sample sizes.
Note that we evaluate the trained models on instances generated with a sample size of 100 regardless of the sample size used during training.
Our results, indicate that there exists a trade-off between different sample sizes. Specifically, our pipeline performs best when training on instances with a sample size of 50, reaching an average gap of $5.15\%$ to the anticipative baseline. Increasing or decreasing the sample size reduces the performance.
We attribute this to the bias-variance trade-off: small instances allow to train a statistical model that is highly accurate on the training set but fails to generalize to the larger instances used during evaluation. Training on large instances on the other hand fails to extract enough structural information to yield an accurate statistical model.

\begin{result}
It is crucial to balance learning accuracy and model generalization when training the ML-CO policy.
\end{result}

\paragraph{Different numbers of training instances}
\begin{table}[ht]
\footnotesize
\centering
\begin{tabular}{c | c c c c c c c c} 
 & \multicolumn{8}{c}{Num. of training instances} \\
 & 1 & 2 & 5 & 10 & 15 & 20 & 25 & 30 \\ 
 \hline\hline
 \textbf{Relative distance to ant. b.:} & 10.19\% & 7.77\% & 7.26\% & 5.89\% & 5.15\% & 5.23\% & 4.50\%  & 4.94\% \\
 \hline
\end{tabular}
\caption{Performance of ML-CO policy when training on different numbers of training instances.}
\label{table:numTrainingInstances}
\end{table}
Table~\ref{table:numTrainingInstances} shows the performance of the ML-CO policy relative to the performance of the anticipative baseline when training the ML-CO policy on different training set sizes.
Our results convey that the performance of our ML-CO policy improves with the training set size. However, the marginal improvement decreases with the training set size. Specifically, performance saturates at a training set size of 10 instances, which shows that our ML-CO policy learns to generalize well even from small training sets. This indicates that our approach requires only a few instances to extract most of the structural information contained in the problem setting, which is a stark contrast to findings from classic supervised learning, which generally requires significantly larger training sets.

\begin{result}
The ML-CO policy only needs few training instances to learn a general policy.
\end{result}

\paragraph{Varying target strategies}
\begin{table}[ht]
\footnotesize
\centering
\begin{tabularx}{\textwidth}{X | c c c c c c c c}
 & \multicolumn{8}{c}{Imitated anticipative strategies} \\
  & best seed & 3600 sec & 900 sec & 300 sec & 240 sec & 180 sec & 120 sec & 60 sec \\ 
 \hline\hline
 \textbf{Relative distance to ant. b.:} & 6.18\% & 5.15\% & 4.79\% & 7.61\% & 6.01\% & 5.25\% & 5.33\% & 5.81\% \\
  \textbf{Average objective value [$\times$10000]:} & 20.49 & 20.56 & 20.70 & 20.78 & 20.80 & 20.87 & 20.97 & 21.21 \\
 \hline
\end{tabularx}
\caption{Performance of ML-CO policy for different target strategies.}
\label{table:imitation}
\end{table}
To investigate the impact of the quality of the solutions imitated by our ML-CO policy, we vary the time limit allocated to the \gls*{HGS} used to derive the underlying anticipative solutions. Here the intuition is as follows: solutions derived with a low time limit should be of lower quality than solutions derived with a high time limit. As the \gls*{HGS} algorithm is subject to random decisions, we further include a training set derived from the best solutions found across 10 runs with a time limit of 3600 seconds each.
Table \ref{table:imitation} shows the performance of the ML-CO policy relative to the performance of the anticipative baseline when training the ML-CO policy on different target solutions with different solution qualities. 
Surprisingly, there is no clear trend between the solution quality of the respective training set and the performance of the ML-CO policy.
This leads to the assumption that an improved quality of the anticipative target solution does not increase the performance of the trained ML-CO policy and therefore the anticipative strategy might not be the best policy to imitate.
\begin{result}
The anticipative strategy might not be the best policy to imitate and there might exist a target solution which yields better performances.
\end{result}

\paragraph{Different statistical models}
\begingroup
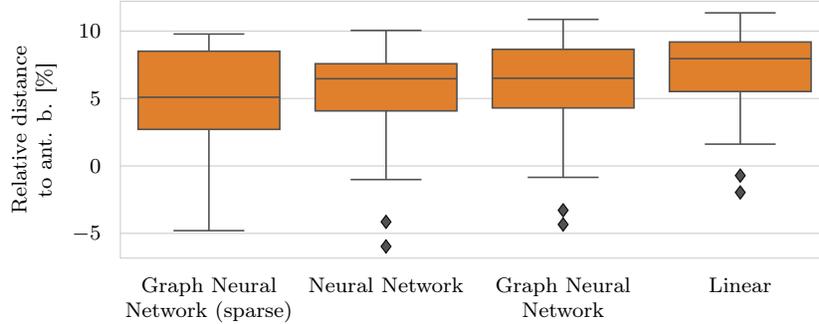
\begin{figure}[ht]
    \centering
    \footnotesize
\begin{tikzpicture}

\definecolor{darkslategray38}{RGB}{38,38,38}
\definecolor{darkslategray81}{RGB}{81,81,81}
\definecolor{lightgray204}{RGB}{204,204,204}
\definecolor{peru22412844}{RGB}{224,128,44}

\begin{axis}[
axis line style={lightgray204},
height=5cm,
scaled y ticks=true,
tick align=outside,
width=11cm,
x grid style={lightgray204},
xmajorticks=false,
xmajorticks=true,
xmin=-0.5, xmax=3.5,
xtick style={color=darkslategray38},
xtick style={draw=none},
xtick={0,1,2,3},
xticklabel style={align=center},
xticklabels={
  Graph Neural \\ Network (sparse),
  Neural Network,
  Graph Neural \\ Network,
  Linear
},
y grid style={lightgray204},
ylabel style={align=center},
ylabel={Relative distance \\ to ant. b. [\%]},
ymajorgrids,
ymajorticks=false,
ymajorticks=true,
ymin=-6.83572755905614, ymax=12.2154870843647,
ytick style={color=darkslategray38},
ytick style={draw=none}
]
\path [draw=darkslategray81, fill=peru22412844, semithick]
(axis cs:-0.4,2.70964450412464)
--(axis cs:0.4,2.70964450412464)
--(axis cs:0.4,8.50606371969132)
--(axis cs:-0.4,8.50606371969132)
--(axis cs:-0.4,2.70964450412464)
--cycle;
\path [draw=darkslategray81, fill=peru22412844, semithick]
(axis cs:0.6,4.07950156378079)
--(axis cs:1.4,4.07950156378079)
--(axis cs:1.4,7.57757278779952)
--(axis cs:0.6,7.57757278779952)
--(axis cs:0.6,4.07950156378079)
--cycle;
\path [draw=darkslategray81, fill=peru22412844, semithick]
(axis cs:1.6,4.2989567660543)
--(axis cs:2.4,4.2989567660543)
--(axis cs:2.4,8.64568853192661)
--(axis cs:1.6,8.64568853192661)
--(axis cs:1.6,4.2989567660543)
--cycle;
\path [draw=darkslategray81, fill=peru22412844, semithick]
(axis cs:2.6,5.51688485962361)
--(axis cs:3.4,5.51688485962361)
--(axis cs:3.4,9.19014499096704)
--(axis cs:2.6,9.19014499096704)
--(axis cs:2.6,5.51688485962361)
--cycle;
\addplot [semithick, darkslategray81]
table {%
0 2.70964450412464
0 -4.79872319532593
};
\addplot [semithick, darkslategray81]
table {%
0 8.50606371969132
0 9.7793845493087
};
\addplot [semithick, darkslategray81]
table {%
-0.2 -4.79872319532593
0.2 -4.79872319532593
};
\addplot [semithick, darkslategray81]
table {%
-0.2 9.7793845493087
0.2 9.7793845493087
};
\addplot [semithick, darkslategray81]
table {%
1 4.07950156378079
1 -1.01029782551158
};
\addplot [semithick, darkslategray81]
table {%
1 7.57757278779952
1 10.0500268225354
};
\addplot [semithick, darkslategray81]
table {%
0.8 -1.01029782551158
1.2 -1.01029782551158
};
\addplot [semithick, darkslategray81]
table {%
0.8 10.0500268225354
1.2 10.0500268225354
};
\addplot [black, mark=diamond*, mark size=2.5, mark options={solid,fill=darkslategray81}, only marks]
table {%
1 -4.14832007418244
1 -5.96976325708246
};
\addplot [semithick, darkslategray81]
table {%
2 4.2989567660543
2 -0.847071893239715
};
\addplot [semithick, darkslategray81]
table {%
2 8.64568853192661
2 10.8646495747511
};
\addplot [semithick, darkslategray81]
table {%
1.8 -0.847071893239715
2.2 -0.847071893239715
};
\addplot [semithick, darkslategray81]
table {%
1.8 10.8646495747511
2.2 10.8646495747511
};
\addplot [black, mark=diamond*, mark size=2.5, mark options={solid,fill=darkslategray81}, only marks]
table {%
2 -3.28904548914846
2 -4.34256331417935
};
\addplot [semithick, darkslategray81]
table {%
3 5.51688485962361
3 1.61580439355315
};
\addplot [semithick, darkslategray81]
table {%
3 9.19014499096704
3 11.349522782391
};
\addplot [semithick, darkslategray81]
table {%
2.8 1.61580439355315
3.2 1.61580439355315
};
\addplot [semithick, darkslategray81]
table {%
2.8 11.349522782391
3.2 11.349522782391
};
\addplot [black, mark=diamond*, mark size=2.5, mark options={solid,fill=darkslategray81}, only marks]
table {%
3 -0.715862136833477
3 -1.96670496369712
};
\addplot [semithick, darkslategray81]
table {%
-0.4 5.09292541419604
0.4 5.09292541419604
};
\addplot [semithick, darkslategray81]
table {%
0.6 6.4643243249672
1.4 6.4643243249672
};
\addplot [semithick, darkslategray81]
table {%
1.6 6.5034690504779
2.4 6.5034690504779
};
\addplot [semithick, darkslategray81]
table {%
2.6 7.95378171309012
3.4 7.95378171309012
};
\end{axis}

\end{tikzpicture}
    \caption{Performance of ML-CO policy using different statistical models.}
    \label{fig:predictors}
\end{figure}
\endgroup

Figure~\ref{fig:predictors} compares the performance of our ML-CO policy when using different statistical models~$\varphi_w$.
All statistical models rely on a feature mapping $\phi$ that maps a state $x^e$ and a request $i$ of this state to a feature vector $\phi(i,x^e)$ in $\mathbb{R}^{|\phi|}$.
Our \emph{linear model} $\varphi_w(x^e) = (w^\top\phi(i,x))_i$ manages to handle variable size inputs and outputs by applying the same linear model $\phi \mapsto w^\top \phi$ independently to each dimension $i$. Similarly, our \emph{neural network} $\varphi_w(x^e) = \big(g_w(\phi(i,x))\big)_i$  applies an auxiliary neural network $g_w$ independently for each dimension $i$. Finally, we propose two \emph{graph neural networks}. Both graph neural networks consider requests as nodes and the connection between requests as edges. We include a regular and a sparsified graph neural network, the latter contains only those edges which are feasible with respect to request time windows and travel times.
Both graph neural networks receive an input vector~$\phi(i,x)$ for each node~$i$, such that $\varphi_w(x^e) = \big(h_w(\phi(i,x)_{i\in I})\big)$.

The linear model performs worst with an average gap of~$6.85$\% to the anticipative baseline while the sparsified graph neural network performs best with a~$5.04$\% gap to the anticipative baseline.
Comparing the performance of the linear model to the neural networks' performance indicates the importance of a feature generator. In fact, using a simple neural network already lowers the gap to the anticipative baseline to $5.15$\%.
As expected, using a graph neural network that calculates structural features further improves the performance of the ML-CO policy. However, the improvement is rather small in comparison to the performance of the neural network. This suggests that most of the structural information is already included in the structured learning approach. 
\begin{result}
Feature-generating statistical models yield the best performing ML-CO policies.
\end{result}

\section{Conclusion}\label{sec:conclusion}

We presented a novel \gls*{CO}-enriched \gls*{ML} pipeline for a dynamic \gls*{VRP} that was introduced in the EURO Meets NeurIPS Vehicle Routing Competition. Specifically, our work contains several methodological contributions and an extensive computational study. From a methodological perspective, we extend \gls*{ML}-based pipelines to objective functions where the dimension of the predicted objective costs does not match the dimension of the decision variables.
These objective functions amend to other problem settings such that we have made them available in the open source library \texttt{InferOpt.jl}. Moreover, we presented the first pipeline that utilizes a metaheuristic component to solve the \gls*{CO}-layer and showed how to carefully design a metaheuristic that finds heuristic solutions which allow to compute approximate gradients. We showed how to train the \gls*{ML}-layer of this pipeline in a supervised learning fashion, i.e., based on a training set derived from an anticipative strategy. We presented a comprehensive numerical study and show that our policy encoded via the \gls*{CO}-enriched \gls*{ML} pipeline outperforms greedy policies by~13.18\% and even monte-carlo policies, which were granted a longer runtime, by~1.57\% in terms of travel time on average. Interestingly, our results point at the fact that counterintuitive to common practice, imitating anticipative strategies can work well for high-dimensional multi-stage stochastic optimization problems, even if the anticipative strategy might not be the best strategy to imitate. 

\subsection*{Acknowledgments}
We thank Wouter Kool and the organizational committee of the EURO Meets NeurIPS 2022 Vehicle Routing Competition for providing the problem setting studied in this paper.

\addcontentsline{toc}{section}{References}
\printbibliography

\end{document}